\documentclass[12pt]{amsart}
\usepackage{latexsym,amsmath,amssymb,graphicx,MnSymbol,color}
\usepackage{tikz}\usetikzlibrary{matrix,arrows,calc}
\graphicspath{{.}{./figs/}}

\theoremstyle{plain}
\newtheorem{thm}{Theorem}[section]
\newtheorem{main}{Theorem}
\newtheorem{lem}[thm]{Lemma}
\newtheorem{prop}[thm]{Proposition}
\newtheorem{cor}[thm]{Corollary}

\theoremstyle{definition}
\newtheorem{defn}[thm]{Definition}
\newtheorem{exmp}[thm]{Example}
\newtheorem{rem}[thm]{Remark}
\newtheorem{quest}[thm]{Question}

\newcommand{\wt}{\widetilde}
\newcommand{\into}{\hookrightarrow}
\newcommand{\onto}{\twoheadrightarrow}
\newcommand{\bull}{\ensuremath{\bullet\ }}
\newcommand{\bc}{\begin{center}}\newcommand{\ec}{\end{center}}
\newcommand{\bt}{\begin{tabular}}\newcommand{\et}{\end{tabular}}

\newcommand{\Z}{\mathbb{Z}}

\newcommand{\R}{\mathbb{R}}

\newcommand{\sph}{\mathbf{S}}

\newcommand{\isom}{\textsc{Isom}}
\newcommand{\lin}{\textsc{Lin}}
\newcommand{\aff}{\textsc{Aff}}
\newcommand{\dir}{\textsc{Dir}}
\newcommand{\mov}{\textsc{Mov}}
\newcommand{\mn}{\textsc{Min}}
\newcommand{\fix}{\textsc{Fix}}
\newcommand{\inv}{\textsc{inv}}

\newcommand{\sym}{\textsc{Sym}}

\newcommand{\cox}{\textsc{Cox}}
\newcommand{\art}{\textsc{Art}}
\newcommand{\midd}{\textsc{Mid}}
\newcommand{\gar}{\textsc{Gar}}
\newcommand{\cryst}{\textsc{Cryst}}
\newcommand{\dart}{\textsc{Art${}^*$}}

\newcommand{\ggar}{G}
\newcommand{\gdiag}{D_w}
\newcommand{\gart}{A}
\newcommand{\gfac}{F_w}
\newcommand{\ghor}{H_w}
\newcommand{\ccryst}{C}
\newcommand{\cdiag}{D}
\newcommand{\ccox}{W}
\newcommand{\cfac}{F}
\newcommand{\chor}{H}

\newcommand{\McLattGeoEquiv}{Proposition~7.5}
\newcommand{\McLattRefl}{Theorem~9.6}

\newcommand{\drawEdge}[2]{\draw[-,thick] #1--#2;}
\newcommand{\drawDashEdge}[2]{\draw[-,dashed,thick] #1--#2;}
\newcommand{\drawArrow}[3]{
  \draw[-,thick] ($#1!.45!#2$)--($#1!.55!#2+#3$);
  \draw[-,thick] ($#1!.45!#2$)--($#1!.55!#2-#3$);
}
\newcommand{\drawTripleEdge}[3]{
  \draw[-,thick] #1--#2; 
  \draw[-,thick,yshift=.5mm] #1--#2; 
  \draw[-,thick,yshift=-.5mm] #1--#2;
  \drawArrow{#1}{#2}{#3}
}
\newcommand{\drawDoubleEdge}[3]{
  \draw[-,thick,yshift=.3mm] #1--#2; 
  \draw[-,thick,yshift=-.3mm] #1--#2;
  \drawArrow{#1}{#2}{#3}
}
\newcommand{\drawDashDoubleEdge}[3]{
  \draw[-,dashed,thick,yshift=.3mm] #1--#2; 
  \draw[-,dashed,thick,yshift=-.3mm] #1--#2;
  \drawArrow{#1}{#2}{#3}
}
\newcommand{\drawInfEdge}[2]{
  \drawDashEdge{#1}{#2} 
  \node[anchor=south] at ($#1!.5!#2$) {$\infty$};
}
\newcommand{\drawDotEdge}[2]{
  \fill ($#1!.3!#2$) circle (.4mm);
  \fill ($#1!.5!#2$) circle (.4mm);
  \fill ($#1!.7!#2$) circle (.4mm);
}
\newcommand{\drawRegDot}[1]{\fill #1 circle (.7mm);\fill[color=black!80] #1 circle (.5mm);}
\newcommand{\drawSpeDot}[1]{\fill #1 circle (1mm);\fill[color=white] #1 circle (.7mm);}
\newcommand{\drawESpDot}[1]{\fill #1 circle (1mm);\fill[color=red!70] #1 circle (.7mm);}

\begin{document}

\title{The structure of euclidean Artin groups}
\author{Jon McCammond}
\address{Dept. of Math., University of California, Santa Barbara, CA 93106} 
\email{jon.mccammond@math.ucsb.edu}
\date{\today}

\begin{abstract}
  The Coxeter groups that act geometrically on euclidean space have
  long been classified and presentations for the irreducible ones are
  encoded in the well-known extended Dynkin diagrams.  The
  corresponding Artin groups are called euclidean Artin groups and,
  despite what one might naively expect, most of them have remained
  fundamentally mysterious for more than forty years.  Recently, my
  coauthors and I have resolved several long-standing conjectures
  about these groups, proving for the first time that every
  irreducible euclidean Artin group is a torsion-free centerless group
  with a decidable word problem and a finite-dimensional classifying
  space.  This article surveys our results and the techniques we use
  to prove them.
\end{abstract}


\maketitle

The reflection groups that act geometrically on spheres and euclidean
spaces are all described by presentations of an exceptionally simple
form and general Coxeter groups are defined by analogy.  These
spherical and euclidean Coxeter groups have long been classified and
their presentations are encoded in the well-known Dynkin diagrams and
extended Dynkin diagrams, respectively.  Artin groups are defined by
modified versions of these Coxeter presentations, and they were
initially introduced to describe the fundamental group of a space
constructed from the complement of the hyperplanes in a complexified
version of the reflection arrangement for the corresponding spherical
or euclidean Coxeter group.  The most basic example of a Coxeter group
is the symmetric group and the corresponding Artin group is the braid
group, the fundamental group of a quotient of the complement of a
complex hyperplane arrangement called the braid arrangement.

The spherical Artin groups, that is the Artin groups corresponding to
the Coxeter groups acting geometrically on spheres, have been well
understood ever since Artin groups themselves were introduced by
Pierre Deligne \cite{De72} and by Brieskorn and Saito \cite{BrSa72} in
adjacent articles in the \emph{Inventiones} in 1972.  One might have
expected the euclidean Artin groups to be the next class of Artin
groups whose structure was well-understood, but this was not to be.
Despite the centrality of euclidean Coxeter groups in Coxeter theory
and Lie theory more generally, euclidean Artin groups have remained
fundamentally mysterious, with a few minor exceptions, for the past
forty years.

In this survey, I describe recent significant progress in the study of
these groups.  In particular, my coauthors and I have succeeded in
clarifying the structure of all euclidean Artin groups.  We do this by
showing that each of these groups is isomorphic to a subgroup of a new
class of Garside groups that we believe to be of independent interest.
The results discussed are contained in the following papers:
``Factoring euclidean isometries'' with Noel Brady \cite{BrMc-factor},
``Dual euclidean Artin groups and the failure of the lattice
property'' \cite{Mc-lattice}, and ``Artin groups of euclidean type''
with Robert Sulway \cite{McSu-artin-euclid}.  The first two are
foundational in nature; the third establishes the main results.  The
structure of this survey follows that of the talks I gave in Durham.
The first part corresponds to my first talk and the second part
corresponds to my second talk.

\part{Factoring euclidean isometries}

I begin with a brief sketch of some elementary facts and known results
about Coxeter groups and Artin groups in order to establish a context
for our results.  The discussion then shifts to a seemingly unrelated
topic: the structure of the poset of all minimum length reflection
factorizations of an arbitrary euclidean isometry.  The connection
between these two disparate topics is rather indirect and its
description is postponed until the second part of the article.

\section{Coxeter groups}

Recall that a group is said to act \emph{geometrically} when it acts
properly discontinuously and cocompactly by isometries, and an action
on euclidean space is \emph{irreducible} if there does not exist a
nontrivial orthogonal decomposition of the underlying space so that
the group is a product of subgroups acting on these subspaces.

\begin{defn}[Spherical Coxeter groups]
  The irreducible \emph{spherical Coxeter groups} are those groups
  generated by reflections that act geometrically and irreducibly on a
  sphere in some euclidean space fixing its center.  The
  classification of such groups is classical and their presentations
  are encoded in the well-known Dynkin diagrams.  The \emph{type} of a
  Dynkin diagram is its name in the Cartan-Killing classification and
  it is \emph{crystallographic} or \emph{non-crystalographic}
  depending on whether or not it extends to a euclidean Coxeter group.
  The crystallographic types consist of three infinite families
  ($A_n$, $B_n = C_n$, and $D_n$) and five sporadic examples ($G_2$,
  $F_4$, $E_6$, $E_7$, and $E_8$).  The non-crystallographic types are
  $H_3$, $H_4$ and $I_2(m)$ for $m \neq 3,4,6$.  The subscript is the
  dimension of the euclidean space containing the sphere on which it
  acts.
\end{defn}

\begin{figure}
  \begin{center}
    \begin{tabular}{cc}
      \includegraphics[scale=.55]{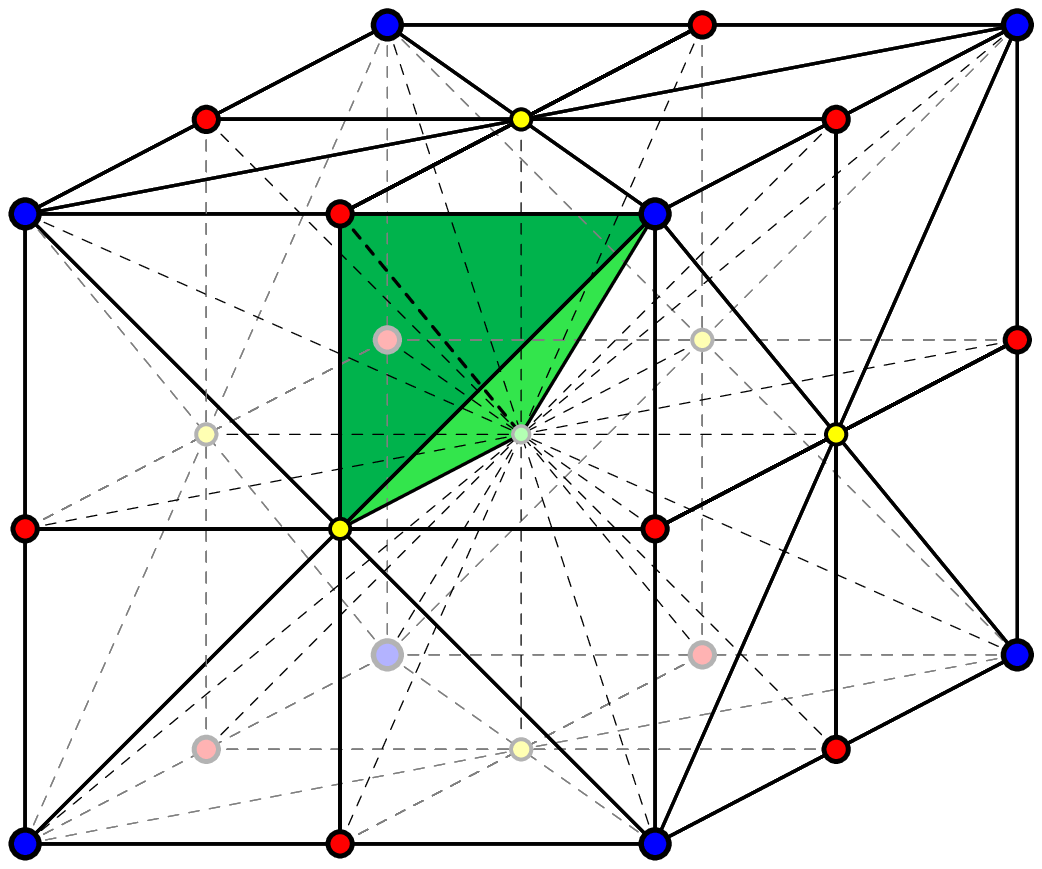} & 
      \includegraphics[scale=.55]{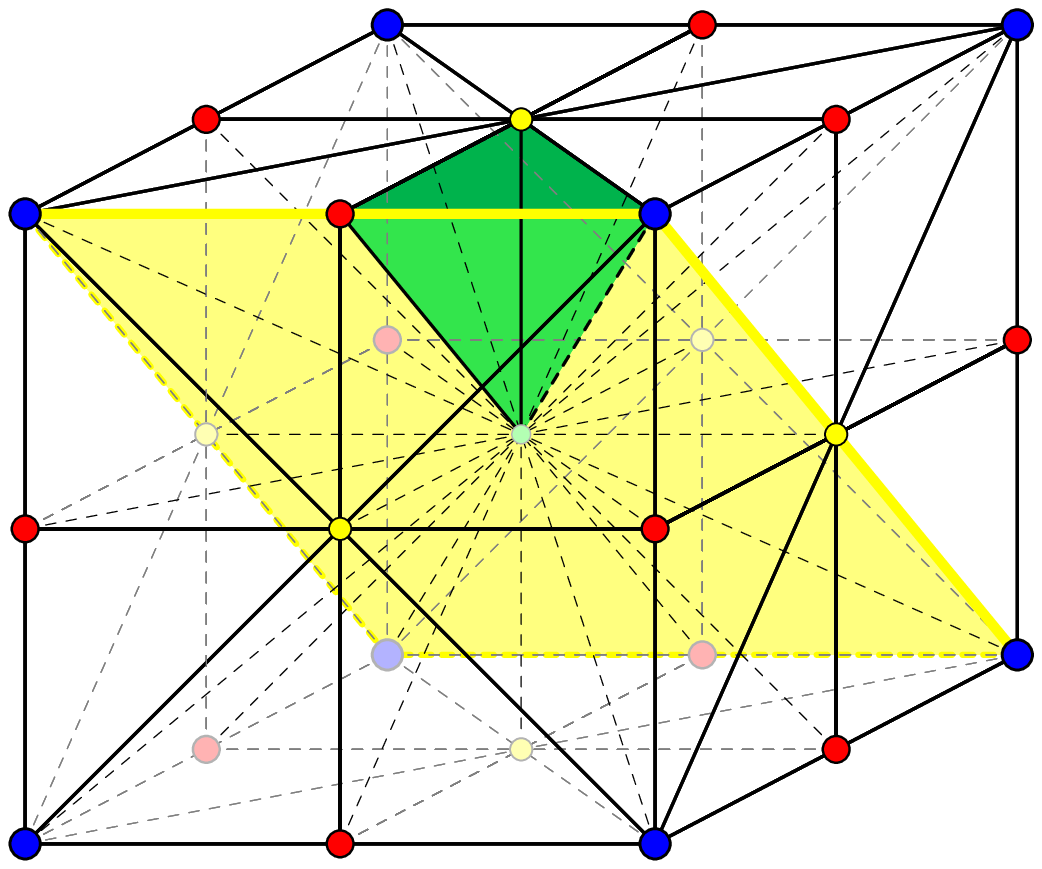} \\
      \includegraphics[scale=.55]{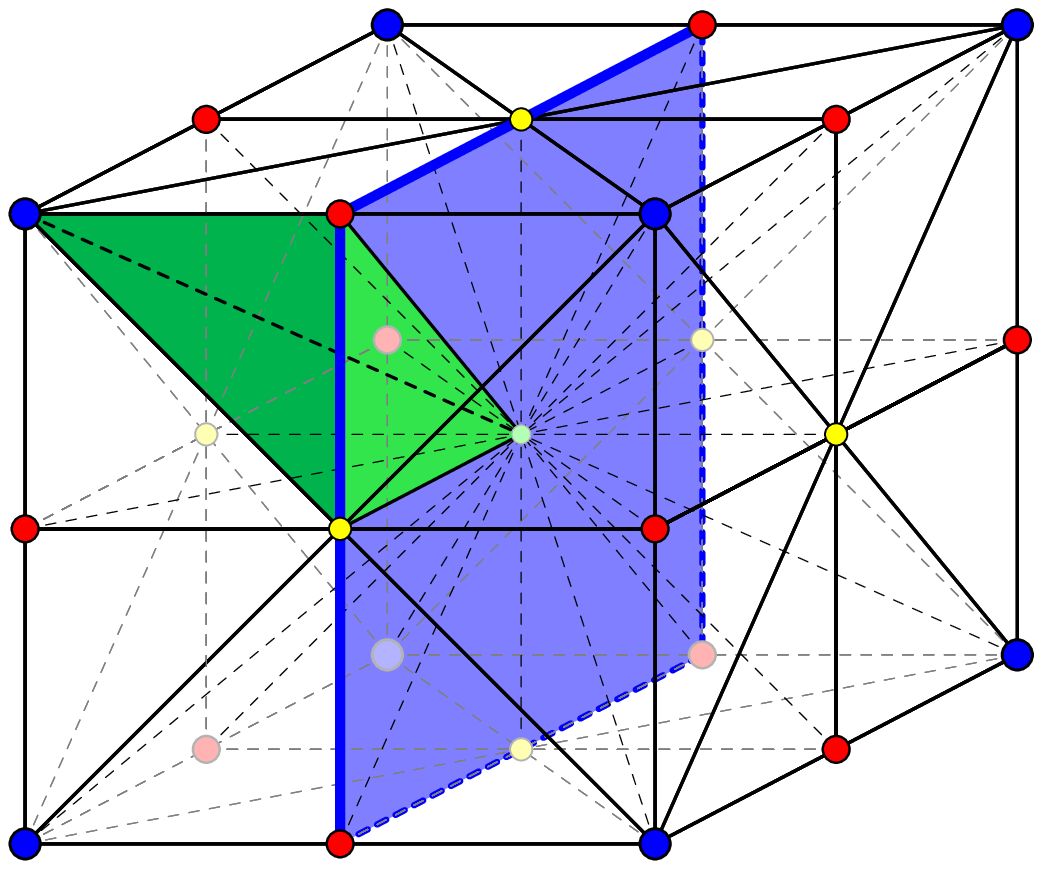} & 
      \includegraphics[scale=.55]{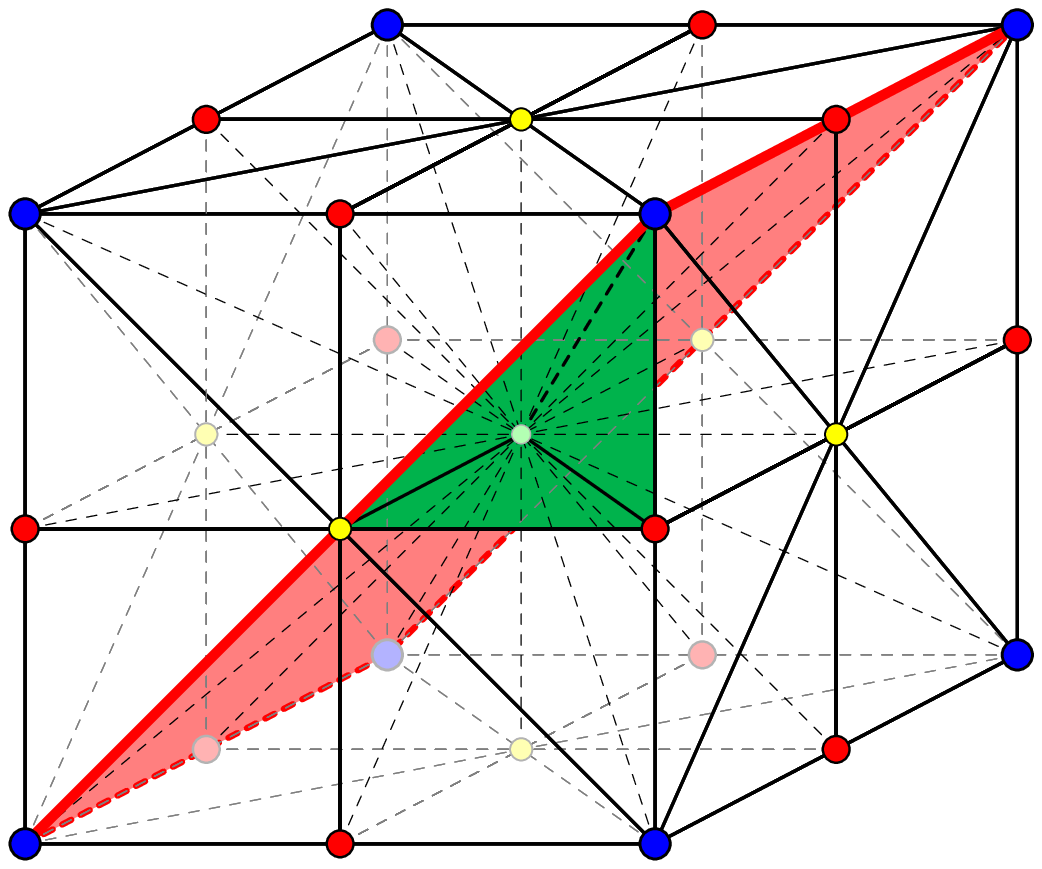} \\
    \end{tabular}
  \end{center}
  \caption{The spherical Coxeter group $\cox(B_3)$.\label{fig:b3}}
\end{figure}

\begin{exmp}[Simplices and cubes]
  The spherical Coxeter groups of types $A$ and $B$ are the best known
  and represent the symmetry groups of regular simplices and
  high-dimensional cubes, respectively.  As groups they are the
  \emph{symmetric groups} and extensions of symmetric groups by
  elementary $2$-groups called a \emph{signed symmetric groups}.  For
  example, the group $\cox(A_3) \cong \sym_4$ is the symmetric group
  of a regular tetrahedron and the group $\cox(B_3) \cong (\Z_2)^3
  \rtimes \sym_3$ and is the group of symmetries of the $3$-cube shown
  in Figure~\ref{fig:b3}.
\end{exmp}

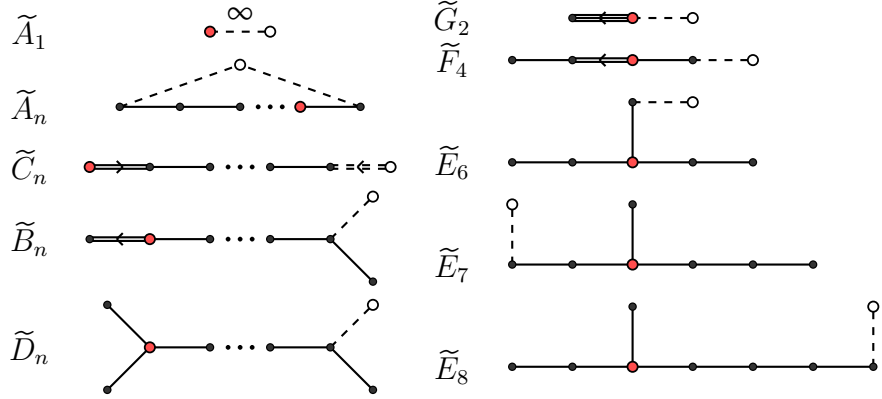
\begin{figure}
  \begin{tabular}{cc}
    \begin{tikzpicture}[scale=.8]
    \begin{scope}[yshift=5.25cm,xshift=2cm]
      \node at (-3,0) {$\wt A_1$};
      \drawInfEdge{(0,0)}{(1,0)}
      \drawSpeDot{(1,0)}
      \drawESpDot{(0,0)}
    \end{scope}
    \begin{scope}[yshift=4cm,xshift=.5cm]
      \node at (-1.5,0) {$\wt A_n$};
      \drawEdge{(0,0)}{(2,0)}
      \drawDashEdge{(0,0)}{(2,.7)}
      \drawDashEdge{(4,0)}{(2,.7)}
      \drawEdge{(3,0)}{(4,0)}
      \drawDotEdge{(2,0)}{(3,0)}
      \foreach \x in {0,1,2,4} {\drawRegDot{(\x,0)}}
      \drawSpeDot{(2,.7)}
      \drawESpDot{(3,0)}
    \end{scope}
    \begin{scope}[yshift=3cm]
      \node at (-1,0) {$\wt C_n$};
      \drawDoubleEdge{(1,0)}{(0,0)}{(0,.1)}
      \drawEdge{(1,0)}{(2,0)}
      \drawDotEdge{(2,0)}{(3,0)}
      \drawEdge{(3,0)}{(4,0)}
      \drawDashDoubleEdge{(4,0)}{(5,0)}{(0,.1)}
      \foreach \x in {1,2,3,4} {\drawRegDot{(\x,0)}}
      \drawESpDot{(0,0)}
      \drawSpeDot{(5,0)}
    \end{scope}
    \begin{scope}[yshift=1.8cm]
      \node at (-1,0) {$\wt B_n$};
      \drawDoubleEdge{(0,0)}{(1,0)}{(0,.1)}
      \drawEdge{(1,0)}{(2,0)}
      \drawDotEdge{(2,0)}{(3,0)}
      \drawEdge{(3,0)}{(4,0)}
      \drawDashEdge{(4,0)}{(4.707,.707)}
      \drawEdge{(4,0)}{(4.707,-.707)}
      \drawRegDot{(4.707,-.707)}
      \foreach \x in {0,2,3,4} {\drawRegDot{(\x,0)}}
      \drawESpDot{(1,0)}
      \drawSpeDot{(4.707,.707)}
    \end{scope}
    \begin{scope}
      \node at (-1,0) {$\wt D_n$};
      \drawEdge{(.293,.707)}{(1,0)}
      \drawEdge{(.293,-.707)}{(1,0)}
      \drawEdge{(1,0)}{(2,0)}
      \drawDotEdge{(2,0)}{(3,0)}
      \drawEdge{(3,0)}{(4,0)}
      \drawDashEdge{(4,0)}{(4.707,.707)}
      \drawEdge{(4,0)}{(4.707,-.707)}
      \drawRegDot{(.293,.707)}
      \drawRegDot{(.293,-.707)}
      \drawRegDot{(4.707,-.707)}
      \foreach \x in {2,3,4} {\drawRegDot{(\x,0)}}
      \drawESpDot{(1,0)}
      \drawSpeDot{(4.707,.707)}
    \end{scope}
  \end{tikzpicture}
  &
  \begin{tikzpicture}[scale=.8]
  \begin{scope}[yshift=5.8cm,xshift=1cm]
    \node at (-2,0) {$\wt G_2$};
    \drawTripleEdge{(0,0)}{(1,0)}{(0,.1)}
    \drawDashEdge{(1,0)}{(2,0)}
    \drawRegDot{(0,0)}
    \drawESpDot{(1,0)}
    \drawSpeDot{(2,0)}
  \end{scope}
  \begin{scope}[yshift=5.1cm]
    \node at (-1,0) {$\wt F_4$};
    \drawDoubleEdge{(1,0)}{(2,0)}{(0,.1)}
    \drawEdge{(0,0)}{(1,0)}
    \drawEdge{(2,0)}{(3,0)}
    \drawDashEdge{(3,0)}{(4,0)}
    \foreach \x in {0,1,3} {\drawRegDot{(\x,0)}}
    \drawESpDot{(2,0)}
    \drawSpeDot{(4,0)}
  \end{scope}
  \begin{scope}[yshift=3.4cm]
    \node at (-1,0) {$\wt E_6$};
    \drawEdge{(0,0)}{(4,0)}
    \drawEdge{(2,0)}{(2,1)}
    \drawDashEdge{(2,1)}{(3,1)}
    \foreach \x in {0,1,3,4} {\drawRegDot{(\x,0)}}
    \drawRegDot{(2,1)}
    \drawESpDot{(2,0)}
    \drawSpeDot{(3,1)}
  \end{scope}
  \begin{scope}[yshift=1.7cm]
    \node at (-1,0) {$\wt E_7$};
    \drawEdge{(0,0)}{(5,0)}
    \drawEdge{(2,0)}{(2,1)}
    \drawDashEdge{(0,0)}{(0,1)}
    \foreach \x in {0,1,3,4,5} {\drawRegDot{(\x,0)}}
    \drawRegDot{(2,1)}
    \drawESpDot{(2,0)}
    \drawSpeDot{(0,1)}
  \end{scope}
  \begin{scope}
    \node at (-1,0) {$\wt E_8$};
    \drawEdge{(0,0)}{(6,0)}
    \drawDashEdge{(6,0)}{(6,1)}
    \drawEdge{(2,0)}{(2,1)}
    \foreach \x in {0,1,3,4,5,6} {\drawRegDot{(\x,0)}}
    \drawRegDot{(2,1)}
    \drawESpDot{(2,0)}
    \drawSpeDot{(6,1)}
  \end{scope}
  \end{tikzpicture}
  \end{tabular}
  \caption{Four infinite families and five sporadic examples.\label{fig:dynkin}}
\end{figure}

\begin{defn}[Euclidean Coxeter groups]
  The irreducible \emph{euclidean Coxeter groups} are the groups
  generated by reflections that act geometrically and irreducibly on
  euclidean space.  The classification of such groups is also
  classical and their presentations are encoded in the extended Dynkin
  diagrams shown in Figure~\ref{fig:dynkin}.  There are four infinite
  families ($\wt A_n$, $\wt B_n$, $\wt C_n$ and $\wt D_n$) and and
  five sporadic examples ($\wt G_2$, $\wt F_4$, $\wt E_6$, $\wt E_7$,
  and $\wt E_6$).  The subscript is the dimension of the euclidean
  space on which it acts.  Removing the white dot and the attached
  dashed edge or edges from the extended Dynkin diagram $\wt X_n$ produces
  the corresponding Dynkin diagram $X_n$.
\end{defn}

These extended Dynkin diagrams index many different objects including
the Artin groups that are our primary focus, but in the present
context, it is more relevant that they index euclidean simplices with
restricted dihedral angles.

\begin{defn}[Euclidean Coxeter simplices]
  Every extended Dynkin diagram encodes a simplex in euclidean space,
  unique up to rescaling, with the following properties: the vertices
  of the diagram are in bijection with the facets of the simplex, i.e.
  its codimension one faces, and vertices $s$ and $t$ in the diagram
  are connected with $0$, $1$, $2$, or $3$ edges iff the corresponding
  facets intersect with a dihedral angle of $\frac{\pi}{2}$,
  $\frac{\pi}{3}$, $\frac{\pi}{4}$, or $\frac{\pi}{6}$, respectively.
  These conventions are sufficient to describe the simplices
  associated to each diagram with one exception: the diagram $\wt A_1$
  corresponds to a $1$-simplex in $\R^1$ whose facets are its
  endpoints.  These do not intersect and this is indicated by the
  infinity label on its unique edge.  The extended Dynkin diagrams
  form a complete list of those euclidean simplices where every
  dihedral angle is of the form $\frac{\pi}{m}$ for some integer $m
  >1$.  We call these \emph{euclidean Coxeter simplices}.
\end{defn}

From these euclidean Coxeter simplices we can recover the
corresponding euclidean Coxeter groups and an associated euclidean
tiling.

\begin{defn}[Euclidean tilings]
  Let $\wt X_n$ be an extended Dynkin diagram and let $\sigma$ be the
  corresponding euclidean $n$-simplex described above.  The group
  generated by the collection of $n+1$ reflections which fix some
  facet of $\sigma$ is the corresponding euclidean Coxeter group $W =
  \cox(\wt X_n)$ and the images of $\sigma$ under the action of $W$
  group tile euclidean $n$-space.  As an illustration, consider the
  extended Dynkin diagram $\wt G_2$.  It represents a euclidean
  triangle with dihedral angles $\frac{\pi}{3}$, $\frac{\pi}{6}$ and
  $\frac{\pi}{2}$ and the euclidean Coxeter group $\cox(\wt G_2)$
  generated by the reflections in its sides is associated with the
  tiling of $\R^2$ by congruent 30-60-90 triangles shown in
  Figure~\ref{fig:g2}.
\end{defn}

\begin{rem}[Spherical analogues]
  For an ordinary Dynkin diagram of type $X_n$, one constructs
  spherical simplex $\sigma$ with similarly restricted dihedral angles
  and recovers the spherical Coxeter group $\cox(X_n)$ as the group
  generated by the reflections in the facets of $\sigma$.  The images
  of $\sigma$ under this action yield a spherical tiling.  This is
  illustrated in Figure~\ref{fig:b3} if one intersects the cell
  structure shown with a small sphere around the center of the cube.
  The cube in the upper left shades a tetrahedron which intersects
  with the small sphere to produce a spherical triangle with dihedral
  angles $\frac{\pi}{3}$, $\frac{\pi}{4}$ and $\frac{\pi}{2}$.  The
  other three cubes illustrate its image under the action of the three
  reflections in its sides.
\end{rem}

And finally a short remark about how the spherical and euclidean cases
relate to the general theory.

\begin{rem}[General Coxeter groups]
  The general theory of Coxeter groups was pioneered by Jacques Tits
  in the early 1960s and the spherical and euclidean Coxeter groups
  are key examples that motivate their introduction.  Coxeter groups
  are defined by simple presentations and in that first unpublished
  paper, Tits proved that every Coxeter group has a faithful linear
  representation preserving a symmetric bilinear form and thus has a
  solvable word problem.  Irreducible Coxeter groups can be coarsely
  classified by the signature of the symmetric bilinear forms they
  preserve and the irreducible spherical and euclidean groups are
  those that preserve positive definite and positive semi-definite
  forms, respectively.
\end{rem}

\section{Artin groups}
As mentioned in the introduction, Artin groups first appear in print
in 1972 in a pair of articles by Pierre Deligne \cite{De72} and by
Brieskorn and Saito \cite{BrSa72}.  Both articles focus on spherical
Artin groups as fundamental groups of spaces constructed from
complements of complex hyperplane arrangements and successfully
analyze their struture using different techniques.  The resulting
presentations resemble Artin's standard presentation for the braid
groups, which is, of course, the most prominent example of a spherical
Artin group.  In the spherical and euclidean context these
presentations are extremely easy to describe.

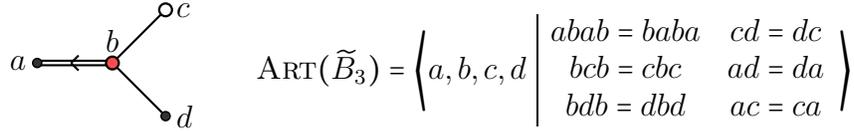
\begin{figure}
  \begin{tabular}{cc}
    \bt{c}
    \begin{tikzpicture}
      \begin{scope}
        \drawDoubleEdge{(0,0)}{(1,0)}{(0,.1)}
        \drawEdge{(1,0)}{(1.707,.707)}
        \drawEdge{(1,0)}{(1.707,-.707)}
        \drawRegDot{(1.707,-.707)}
        \foreach \x in {0} {\drawRegDot{(\x,0)}}
        \drawESpDot{(1,0)}
        \drawSpeDot{(1.707,.707)}
        \node at (0,0) [anchor=east] {$a$};
        \node at (1,0) [anchor=south] {$b$};
        \node at (1.707,.707) [anchor=west] {$c$};
        \node at (1.707,-.707) [anchor=west] {$d$};
      \end{scope}
    \end{tikzpicture}
    \et
    &
    \bt{c} $\art(\wt B_3) = \left< a,b,c,d\  \begin{array}{|cc} abab = baba & cd=dc\\ 
      bcb=cbc & ad=da\\ bdb=dbd & ac=ca \\ \end{array} \right>$ \et
  \end{tabular}
  \caption{The $\wt B_3$ diagram and the presentation for the
    corresponding euclidean Artin group.\label{fig:artb3}}
\end{figure}

\begin{defn}[Euclidean Artin groups]
  Let $\wt X_n$ be an extended Dynkin diagram.  The standard
  presentation for the Artin group of type $\wt X_n$ has a generator
  for each vertex and at most one relation for each pair of vertices.
  More precisely, if $s$ and $t$ are vertices connected by $0$, $1$,
  $2$, or $3$ edges, then the presentation contains the relation
  $st=ts$, $sts=tst$, $stst=tsts$ or $ststst=tststs$ respectively.
  And finally, in the case of $\wt A_1$, the edge labeled $\infty$
  indicates that there is no relation corresponding to this pair of
  vertices.  As an illustration, Figure~\ref{fig:artb3} shows the
  extended Dynkin diagram of type $\wt B_3$ along with the explicit
  presentation for the corresponding euclidean Artin group $\art(\wt
  B_3)$.
\end{defn}

General Artin groups are defined by similarly simple presentations
encoded in the same diagrams as general Coxeter groups and then
coarsely classified in the same way.  Given the centrality of
euclidean Coxeter groups and the elegance of their structure, one
might have expected euclidean Artin groups to be well understood
shortly thereafter.  It is now 40 years later and these groups are
still revealing their secrets.

\begin{defn}[Four conjectures]
  In a recent survey article, Eddy Godelle and Luis Paris highlight
  how little we know about general Artin groups by highlighting four
  basic conjectures that remain open \cite{GoPa-basic}.  Their four
  conjectures are:
  \begin{enumerate}
    \item[(A)] All Artin groups are torsion-free.
    \item[(B)] Every non-spherical irreducible Artin group has a trivial center.
    \item[(C)] Every Artin group has a solvable word problem.
    \item[(D)] All Artin groups satisfy the $K(\pi,1)$ conjecture.
  \end{enumerate}
  Godelle and Paris also remark that these conjectures remain open and
  are a ``challenging question'' even in the case of the euclidean
  Artin groups.  These are precisely the conjectures that my
  collaborators and I set out to resolve.
\end{defn}

There are a few euclidean Artin groups with a well-undestood
structure.  The earliest results are by Craig Squier.

\begin{exmp}[planar Artin groups]
  In a 1987 article, Squier successfully analyzes the structure of the
  three irreducible euclidean Artin groups $\art(\wt A_2)$, $\art(\wt
  C_2)$, and $\art(\wt G_2)$ that correspond to the three irreducible
  euclidean Coxeter groups which act geometrically on the euclidean
  plane \cite{Squier87}.  He works directly with the presentations and
  analyzes them as amalgamated products and HNN extensions of
  well-known groups.  This technique does not appear to generalize to
  other euclidean Artin groups.
\end{exmp}

\begin{figure}
  \begin{center}
    \includegraphics[scale=.8]{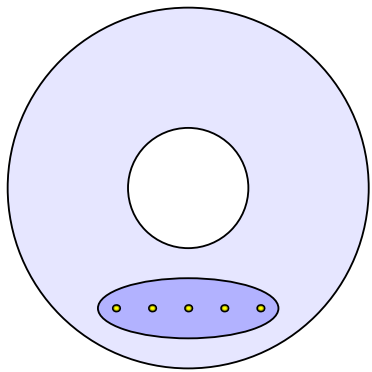}
  \end{center}
  \caption{Five punctures in a disk inside an
    annulus.\label{fig:annular-braids}}
\end{figure}

A second result is the consequence of an unusual embedding of a
euclidean Artin group into a spherical Artin group.

\begin{exmp}[Annular braids]
  It has been repeatedly observed that the euclidean Artin group
  $\art(\wt A_n)$ embeds into the spherical Artin group
  $\art(B_{n+1})$ and is, in fact, part of a short exact sequence
  \[\art(\wt A_n) \into \art(B_{n+1}) \onto \Z\] which greatly clarifies
  its structure \cite{tD98,Al02,KePe02,ChPe03}.  The group $\art(B_n)$
  is sometimes called the \emph{annular braid group} because it can be
  interpreted as the braid group of the annulus \cite{Bi74}.  If one
  selects a disk in the annulus containing all the punctures, as shown
  in Figure~\ref{fig:annular-braids}, then the path traced by each
  puncture, viewed as a path that starts and ends in the disk, has a
  winding number.  The sum of these individual winding numbers is a
  global winding number for each element of $\art(B_{n+1})$, and this
  assignment of a global winding number is a group homomorphism onto
  $\Z$ with $\art(\wt A_n)$ as its kernel.  In other words, the group
  $\art(\wt A_n)$ is the subgroup of annular braids with global
  winding number $0$.
\end{exmp}

And finally, there are two recent results due to Fran\c{c}ois Digne. 

\begin{exmp}[Garside structures]
  Digne showed that the groups $\art(\wt A_n)$ and $\art(\wt C_n)$
  have infinite-type Garside structures \cite{Digne06, Digne12}.  In
  the first article Digne uses the embedding $\art(\wt A_n) \into
  \art(B_{n+1})$ to show that the euclidean Artin groups of type $\wt
  A_n$ have infinite-type Garside structures and in the second he uses
  a delicate analysis of the some maps relating type $C$ and type $A$
  to show that the euclidean Artin groups of type $\wt C_n$ also has
  an infinite-type Garside structure.  Our approach to arbitrary
  euclidean Artin groups is closely related to Digne's work and the
  second part of the article contains a more detailed description of
  Garside structures and their uses.
\end{exmp}

To my knowledge, these euclidean Artin groups, i.e. the ones of type
$\wt A_n$, $\wt C_n$, and $\wt G_2$, are the only ones whose structure
was previously fully understood.  In fact, one of the main
frustrations in the area is the stark contrast between the utter
simplicity of the presentations involved and the fact that we
typically know very little about the groups they define.  

For example, all four conjectures identified by Godelle and Paris were
open for the group $\art(\wt B_3)$ shown in Figure~\ref{fig:artb3} --
including a solution to its word problem -- until 2010 when my
Ph.D. student Robert Sulway analyzed its structure as part of his
dissertation \cite{Sulway-diss}.  As an extension of Sulway's work, he
and I are now able to give positive solutions to Conjectures (A), (B)
and (C) for all euclidean Artin groups and we also make some progress
on Conjecture (D).  We prove, in particular, that every irreducible
euclidean Artin group $\art(\widetilde X_n)$ is a torsion-free
centerless group with a solvable word problem and a finite-dimensional
classifying space.  Our proofs rely heavily on the structure of
intervals in euclidean Coxeter groups and other euclidean groups
generated by reflections, and so we now shift our attention to
structural aspects of the set of all factorizations of a euclidean
isometry into reflections.

\section{Isometries}

Every euclidean isometry can be built out of reflections and the
Cayley graph of the euclidean isometry group with respect to this
natural reflection generating set has bounded diammeter.  This follows
from a fact that most mathematicians learn early on in their
education: every isometry of $n$-dimensional euclidean space is a
product of at most $n+1$ reflections.  The goal of the next few
sections is to describe in some detail the structure of the portion of
this Cayley graph between the identity and a fixed euclidean isometry.
We begin with a coarse classification of euclidean isometries and
their basic invariants following the approach taken in
\cite{BrMc-factor}.  The first step is elementary but important for
conceptual clarity: we make a sharp distinction between points and
vectors.

\begin{figure}
  \bc\bt{ccc}
  \bt{c}\includegraphics[scale=.7]{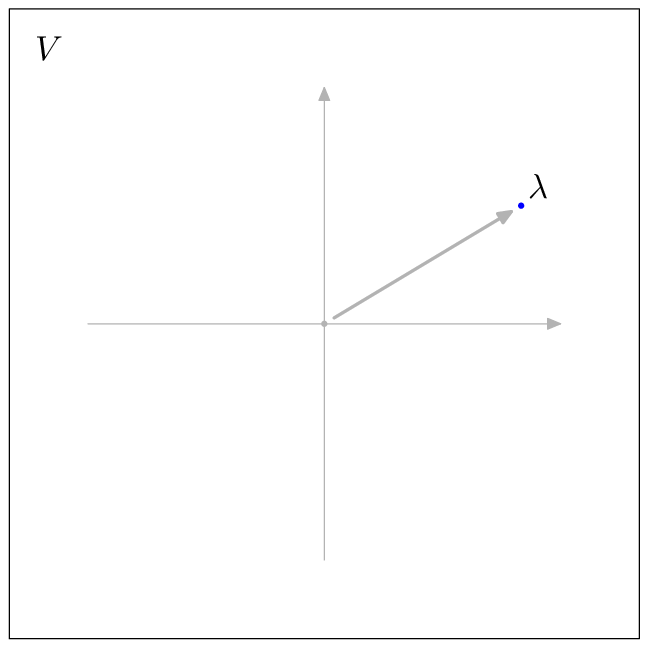}\et 
  & $\curvearrowright$ &
  \bt{c}\includegraphics[scale=.7]{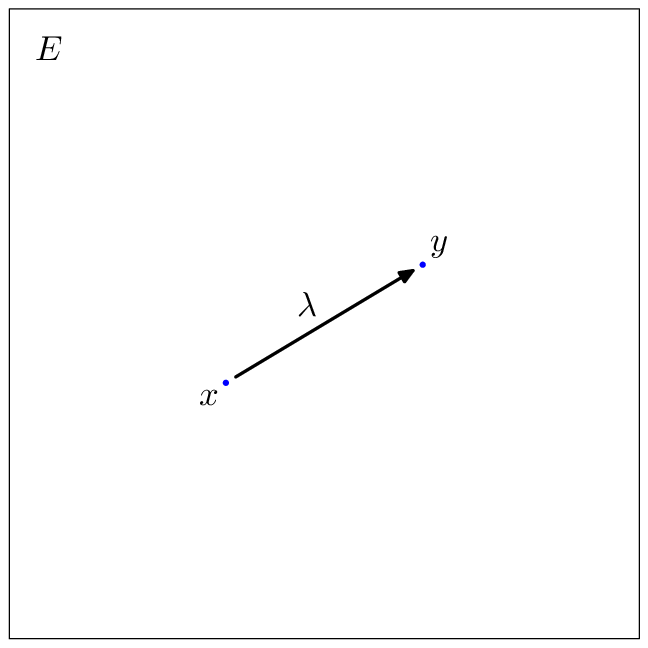}\et 
  \et\ec
  \caption{Vectors acting on points.\label{fig:vec-pt}}
\end{figure}

\begin{defn}[Points and vectors]
  Let $V$ be a vector space with a simple transitive action on a set
  $E$ as shown in Figure~\ref{fig:vec-pt}.  We call $E$ an
  \emph{affine space}, the elements of $E$ are called \emph{points}
  and the elements of $V$ are called \emph{vectors}.  The main
  difference between $V$ and $E$ is that $E$ is essentially a vector
  space with no distinguished point identified as its origin.
\end{defn}

Both $V$ and $E$ have a natural collections of subspaces which are
used to defines the basic invariants of euclidean isometries.

\begin{defn}[Subspaces]
  A subset of $V$ is \emph{linear} if it is closed under linear
  combination.  A subset of $V$ or $E$ is \emph{affine} if for every
  pair of elements in the subset, the line through these elements is
  also in the subset.  Thus the vector space $V$ has \emph{linear
    subspaces} through the origin and other \emph{affine subspaces}
  not through the origin.  The affine space $E$ only has affine
  subspaces.  For any affine subspace $B \subset E$, vectors between
  points in $B$ form a linear subspace $\dir(B) \subset V$ called its
  \emph{space of directions}.
\end{defn}

Posets are obtained by ordering these natural subspaces by inclusion.

\begin{defn}[Poset structure]
  The linear subspaces of $V$ ordered by inclusion define a poset
  $\lin(V)$ which is a graded, bounded, self-dual lattice.  The affine
  subspaces of $E$ ordered by inclusion define a poset $\aff(E)$ which
  is graded and bounded above, but not bounded below, not self-dual
  and not a lattice.  Also note that there is a well-defined
  rank-preserving map $\aff(E) \onto \lin(V)$ that sends $B$ to
  $\dir(B)$.
\end{defn}

If one equips $V$ with a positive definite inner product, this induces
a euclidean metric on $E$ and a corresponding set of euclidean
isometries that preserve this metric.  

\begin{defn}[Basic invariants]
  Let $w$ be an isometry of the euclidean space $E$.  Its
  \emph{move-set} is the subset $\mov(w) \subset V$ of all the motions
  that its points undergo.  In symbols,
  \[\mov(w) = \{ w(x)-x \in V \mid x\in E\}\]
  and it is easy to show that $\mov(w)$ is an affine subspace of $V$.
  As an affine subspace, $\mov(w)$ is a translation of a linear
  subspace.  If $U$ denotes the unique linear subspace of $V$ which
  differs from $\mov(w)$ by a translation and $\mu$ is the unique
  vector in $\mov(w)$ closest to the origin, then we call $U + \mu$
  the \emph{standard form} of $\mov(w)$.  The points in $E$ that
  undergo the motion $\mu$ are a subset $\mn(w) \subset E$ called the
  \emph{min-set} of $w$ and it is also easy to show that $\mn(w)$ is
  an affine subspace of $E$.  We call these the \emph{basic
    invariants} of $w$.
\end{defn}

Euclidean isometries naturally divide into two types.

\begin{defn}[Elliptic and hyperbolic]
  Let $w$ be a euclidean isometry.  If its move-set $\mov(w)$ includes
  the origin, then $\mu$ is trivial, and its min-set $\mn(w)$ is also
  its fix-set $\fix(w)$.  Under these equivalent conditions $w$ is
  called \emph{elliptic}.  Otherwise, $w$ is called \emph{hyperbolic}.
\end{defn}

\begin{figure}
  \bc\bt{ccc}
  \bt{c}\includegraphics[scale=.7]{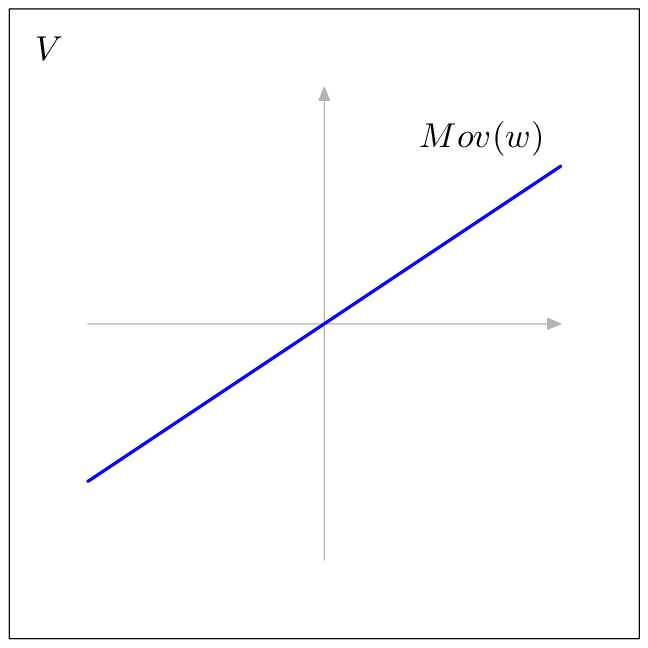}\et 
  & \hspace*{0em} &
  \bt{c}\includegraphics[scale=.7]{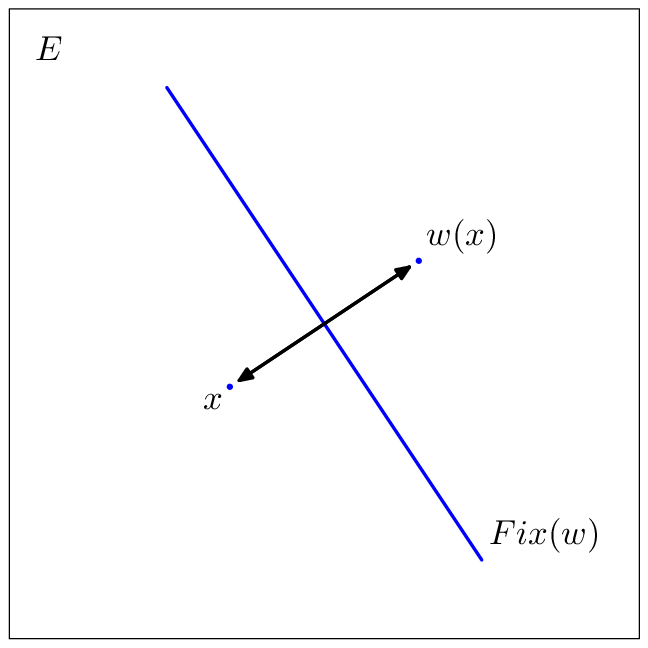}\et 
  \et\ec
  \caption{Basic invariants of a reflection.\label{fig:reflection}}
\end{figure}

The simplest euclidean isometries are reflections and translations.

\begin{defn}[Translations]
  For each vector $\lambda \in V$ there is a \emph{translation}
  isometry $t_\lambda$ whose min-set is all of $E$ and whose move-set
  is the single point $\{\lambda\}$.  So long as $\lambda$ is
  nontrivial, $t_\lambda$ is a hyperbolic isometry.  
\end{defn}

\begin{defn}[Reflections]
  For each \emph{hyperplane} $H$ in $E$ (an affine subspace of
  codimension $1$) there is a unique nontrivial isometry $r$ fixing
  $H$ called a \emph{reflection}.  It is elliptic with fix-set $H$ and
  its move-set is a line through the origin in $V$.  We call any
  nontrivial vector $\alpha$ in this line a \emph{root} of $r$.  The
  basic invariants of a typical reflection are shown in
  Figure~\ref{fig:reflection}.
\end{defn}

\begin{figure}
  \bc\bt{ccc}
  \bt{c}\includegraphics[scale=.7]{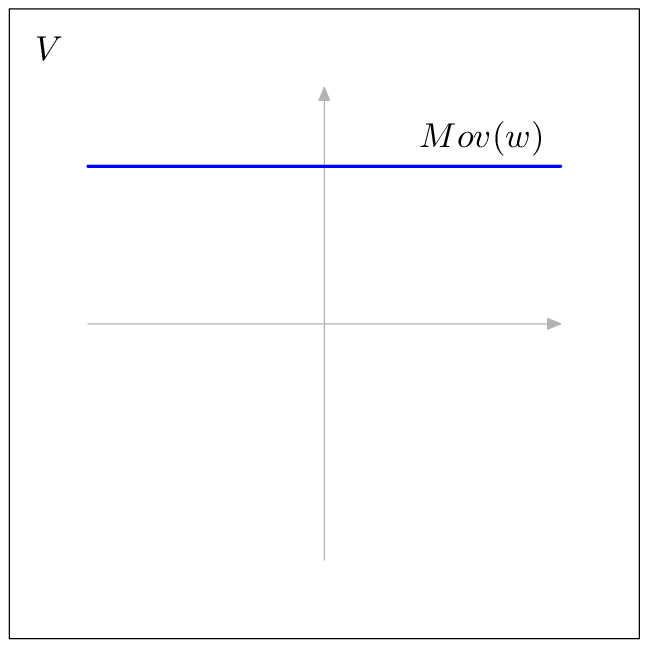}\et 
  & \hspace*{0em} &
  \bt{c}\includegraphics[scale=.7]{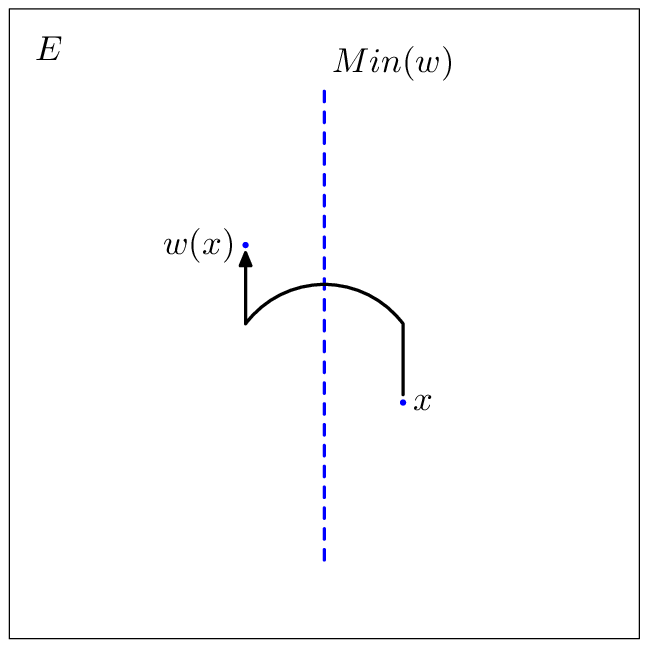}\et 
  \et\ec
  \caption{Basic invariants of a glide reflection.\label{fig:glide}}
\end{figure}

A more interesting example which better illustrates these ideas is
given by a glide reflection.  The move-set of a glide reflection such
as the one shown in Figure~\ref{fig:glide} is a non-linear affine line
in $V$.  It has a unique point $\mu$ closest to the origin and the
points in $E$ which undergo the motion $\mu$ are those on its min-set,
also known as its glide axis.

\section{Intervals}

Let $V$ be an $n$-dimensional vector space over $\R$, let $E$ be an
$n$-dimensional euclidean space on which it acts, and let $L =
\isom(E)$ be the Lie group of euclidean isometries of $E$.  The
structure I want to describe is the portion of the Cayley graph of $L$
generated by its reflections between the identity element and a fixed
euclidean isometry $w$.  We begin by recalling that in any metric
space there is a notion of ``betweenness'' which can used to construct
intervals that are posets.

\begin{defn}[Intervals in metric spaces]
  In any metric space a point $z$ is said to be \emph{between} points
  $x$ and $y$ when the triangle inequality becomes an equality, that
  is when $d(x,z) + d(z,y) = d(x,y)$.  The set of all points between
  $x$ and $y$ form the \emph{interval} $[x,y]$ and this set can be
  given a partial ordering by defining $z \leq w$ if and only if
  $d(x,z) + d(z,w) + d(w,y) = d(x,y)$.
\end{defn}

\begin{figure}
  \begin{center}
    \includegraphics[scale=.5]{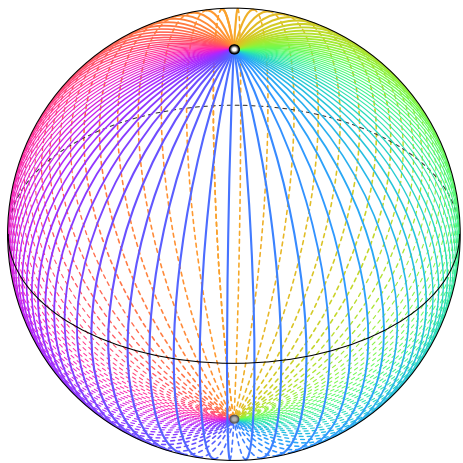}
  \end{center}
  \caption{The interval between antipodal points on a
    $2$-sphere.\label{fig:sphere}}
\end{figure}

As an illustration, consider the unit $2$-sphere with standard angle
metric.  If $x$ and $y$ are not antipodal, the only points between $x$
and $y$ are those along the unique geodesic connecting them with the
usual ordering of an interval in $\R$.  If, however, $x$ and $y$ are
antipodal, say $x$ is the south pole and $y$ is the north pole, then
all points on $\sph^2$ are between $x$ and $y$, the interval $[x,y]$
is all of $\sph^2$ and its ordering is that $z < w$ iff they lie on a
common longitude line with the latitude of $z$ below that of $w$.  See
Figure~\ref{fig:sphere}.  When Cayley graphs are viewed as metric
spaces, they can be used to construct intervals.

\begin{defn}[Intervals in groups]
  Let $G$ be a group with a fixed symmetric generating set $S$.  If we
  assign the elements of $S$ positive weights and the set of all
  possible weights is a discrete subset of the reals, then $G$ can be
  viewed as a metric space where the distance $d(g,h)$ is calculated
  as the minimum total length of a path in the Cayley graph from $v_g$
  to $v_h$.  One convention is to assign every generator a weight of
  $1$ and for finite generating sets the discreteness condition is
  always true, but for infinite generating sets with varying weights,
  some condition is needed so that the infimum of distances between
  two points is achieved by some actual path in the Cayley graph.  Let
  $[g,h]^G$ denote the portion of the Cayley graph between $g$ and
  $h$, by which I mean the union of all the minimal length directed
  paths from $v_g$ to $v_h$.  This is an edge-labeled directed graph
  which also encodes the poset structure.
\end{defn}

Since Cayley graphs are homogeneous, $d(g,h) = d(1,g^{-1}h)$ and the
interval $[g,h]^G$ is isomorphic as an edge-labeled directed graph to
the interval $[1,g^{-1}h]^G$.  Thus it is sufficient to restrict our
attention to distances from the identity and intervals of the form
$[1,g]^G$.  Note that the formulas given above are for right Cayley
graphs with their natural left group action; for left Cayley graphs
$[g,h]^G \cong [1,hg^{-1}]^G$.

\begin{quest}[Euclidean intervals]
  In this language our goal is to describe the poset structure of
  intervals in the Lie group $L = \isom(E)$ of all euclidean
  isometries generated by its full set of reflections with each
  reflection given unit weight.  Questions one might ask include: What
  are the possible poset structures for these intervals $[1,w]^L$? To
  what extent is this poset structure independent of $w$?  Are these
  posets lattices? i.e. do well-defined meets and joins always exist?
  A good test case for these questions is when $w$ is a loxodromic
  ``corkscrew'' isometry of $\R^3$, to which we return at the end of
  Section~\ref{sec:models}.
\end{quest}

In \cite{BrMc-factor} Noel Brady and I answer these questions by
completely characterizing the poset structure of all euclidean
intervals.  Our motivation was to create a technical tool that could
be used to construct dual presentations of euclidean Artin groups, a
process described in Section~\ref{sec:dual-artin}.  The first step is
to understand how far an isometry $w$ is from the identity in this
Cayley graph, i.e. its \emph{reflection length} $\ell_R(w)$, and this
is the content of a classical result known as Scherk's theorem
\cite{Sch50}.

\begin{thm}[Reflection length]
  Let $w$ be a euclidean isometry with a $k$-dimensional move-set.  If
  $w$ is elliptic, its reflection length is $k$.  If $w$ is
  hyperbolic, its reflection length is $k+2$.
\end{thm}

From there we build up an understanding of how the basic invariants of
a euclidean isometry change when it is multiplied by a reflection.
The following lemma is one of the results we establish.

\begin{lem}
  Suppose $w$ is hyperbolic with $\ell_R(w)=k$ and $\mov(w) = U+\mu$
  in standard form, $r$ is a reflection with root $\alpha$ and let
  $U_\alpha$ denote the span of $U \cup \{\alpha\}$.\\ \bull If
  $\alpha \in U$ then $rw$ is hyperbolic with $\ell_R(rw) =
  k-1$.\\ \bull If $\alpha \not \in U$ and $\mu \in U_\alpha$ then
  $rw$ is elliptic with $\ell_R(rw) = k-1$.\\ \bull If $\alpha \not
  \in U$ and $\mu \not \in U_\alpha$ then $rw$ is hyperbolic with
  $\ell_R(rw) = k+1$.
\end{lem}

What I would like the reader to notice is that the geometric
relationships between the basic invariants of $w$ and the basic
invariants of $r$ combine to determine key properties of the basic
invariants of $rw$.  This type of detailed information makes it
possible to prove results such as the following:

\begin{prop}[Elliptic intervals]\label{prop:ell-int}
  Let $w$ be an elliptic isometry with $\mov(w) = U \subset V$.  The
  map $u \mapsto \mov(u)$ creates a poset isomorphism $[1,w]^L \cong
  \lin(U)$.  In particular, $[1,w]^L$ is a lattice.
\end{prop}

Alternatively the map $u \mapsto \fix(u)$ gives a poset isomorphism
with the affine subspaces containing $\fix(w)$ under reverse
inclusion.  Proofs of this result can be found \cite{Sch50},
\cite{BrWa02b} or \cite{BrMc-factor}.  The most remarkable aspect of
this proposition is that the structure of the interval $[1,w]^L$ only
depends on the fact that $w$ is elliptic and the dimension of its
move-set (or equivalently the codimension of its fix-set); it is
otherwise independent of $w$ itself.  In other words, the fix-set of
$w$ completely determines the order structure of the interval.

\section{Models\label{sec:models}}

The main new result established in \cite{BrMc-factor} is an analysis
of the structure of euclidean intervals for hyperbolic isometries.  To
describe these intervals, we first define an abstract poset which
mimics the basic invariants of euclidean isometries.

\begin{defn}[Global Poset]
  Let $E$ be an $n$-dimensional euclidean space and let $V$ be the
  $n$-dimensional vector space that acts on it.  We construct a poset
  $P$ called the \emph{global poset} with two types of elements: it
  has an element we call $h^M$ for each nonlinear affine subspace $M
  \subset V$ and an element we call $e^B$ for each affine subspace $B
  \subset E$.  The ordering of these elements is defined as follows:
  \[ \begin{array}{l} 
    h^M \geq h^{M'} \text{ iff } M \supset M'\\  
    e^B \geq e^{B'} \text{ iff } B \subset B' \\
    h^M > e^B \text{ iff } M^\perp \subset \textsc{Dir}(B)\\ 
    \text{no $e^B$ is ever above $h^M$}
  \end{array}\]
\end{defn}

Next, we define a map from the Lie group $L = \isom(E)$ to the global
poset $P$.

\begin{defn}[Invariant map]
  For each euclidean isometry $w$, the \emph{invariant map} assigns an
  element of $P$ based on its type and its basic invariants.  More
  precisely, the invariant map $\inv:L \to P$ is defined by setting
  $\inv(u) = h^{\mov(u)}$ when $u$ is hyperbolic and $\inv(u) =
  e^{\fix(u)}$ when $u$ is elliptic.
\end{defn}

One reason to introduce the poset $P$ and the map $\inv$ is that there
is a way to use distance from the identity to turn the Lie group $L$
into a poset and under this ordering the invariant map is a
rank-preserving poset map.  It is, however, far from injective as can
be seen from the fact that all rotations which fix the same subspace
are sent to the same element of $P$.  Because $\inv:L\to P$ is a
well-defined map between posets, it sends the elements below $w$ to
the elements below $\inv(w)$.  The former are the intervals $[1,w]^L$.
The latter are what we called model posets.

\begin{defn}[Model posets]
  For each affine subspace $B \subset E$, let $P^B$ denote the poset
  of elements below $e^B$ in the global poset $P$.  Similarly, for
  each nonlinear affine subspace $M \subset V$, let $P^M$ denote the
  poset of elements below $h^M$ in global poset $P$.  We call these
  our \emph{model posets}.  Finally, let $P(w)$ be the model poset of
  elements below $\inv(w)$ in $P$.
\end{defn}

As noted above, the invariant map sends elements in the interval
$[1,w]^L$ is elements in the model poset $P(w)$.  In fact, one of the
main results in \cite{BrMc-factor} is that these restrictions of the
invariant map are poset isomorphisms.

\begin{thm}[Models for euclidean intervals]\label{thm:euc-int}
  For each isometry $w \in L$, the map $u \mapsto \inv(u)$ is a poset
  isomorphism between the interval $[1,w]^L$ and the model poset
  $P(w)$.
\end{thm}

When $w$ is elliptic, this reduces to the previously known
Proposition~\ref{prop:ell-int}, but when $w$ is hyperbolic this result
is new.  Theorem~\ref{thm:euc-int} allows attention to shift away from
the isometries themselves and to focus instead on these model posets
defined purely in terms of the affine subspaces of $V$ and $E$.  In
particular, we are able to understand the structure of euclidean
intervals well enough that we can determine when meets and joins
exist.

\begin{cor}[Lattice failure]
  Let $w\in L$ be a euclidean.  The interval $[1,w]^L$ is not a
  lattice iff $w$ is a hyperbolic isometry and its move-set has
  dimension at least $2$. All other intervals are lattices.
\end{cor}

In \cite{BrMc-factor} we give an explicit characterization of where
these failures occur.  For the applciation to euclidean Artin groups
it is sufficient to describe these failure when $w$ is a hyperbolic
isometry of maximal reflection length.  For such a $w$, its min-set is
a line in $E$ and its move-set is a nonlinear affine hyperplane
(i.e. an affine subspace of codimension~$1$) in $V$.  We call the
direction of its min-set \emph{vertical} and all of the orthogonal
directions \emph{horizontal}.  More generally we call any motion with
a non-trivial vertical component vertical.  One consequence of
Theorem~\ref{thm:euc-int} is that there is exactly one elliptic in
$[1,w]^L$ for each affine subspace $M \subset E$ and exactly one
hyperbolic for each affine subspace of $\mov(w) \subset V$.  Using the
model poset structure as a guide we coarsely partition the elements in
the interval $[1,w]^L$ into a grid with three rows.

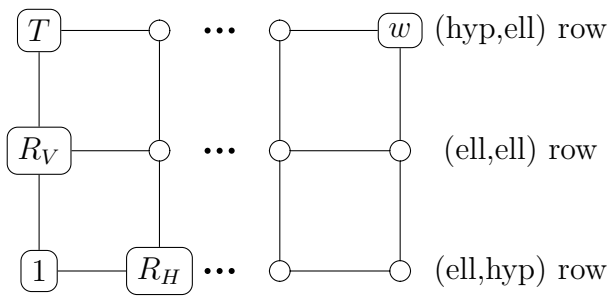
\begin{figure}
  \bc\begin{tikzpicture}[node distance=16mm, auto]
    \begin{scope}
      \tikzstyle{every node}=[rounded corners,draw]
      \node(00){$1$}; 
      \node[right of=00](10){$R_H$}; 
      \node[right of=10](20){}; 
      \node[right of=20](30){}; 

      \node[above of=00](01){$R_V$}; 
      \node[right of=01](11){}; 
      \node[right of=11](21){}; 
      \node[right of=21](31){}; 

      \node[above of=01](02){$T$}; 
      \node[right of=02](12){}; 
      \node[right of=12](22){}; 
      \node[right of=22](32){$w$}; 
    \end{scope}

    \node[right of=30](r1) {(ell,hyp) row};
    \node[right of=31](r2) {(ell,ell) row};
    \node[right of=32](r3) {(hyp,ell) row};

    \draw[-](00)--(10); \draw[-](20)--(30); 
    \draw[-](01)--(11); \draw[-](21)--(31); 
    \draw[-](02)--(12); \draw[-](22)--(32); 

    \fill ($(10)!.4!(20)$) circle (.4mm);
    \fill ($(10)!.5!(20)$) circle (.4mm);
    \fill ($(10)!.6!(20)$) circle (.4mm);
    \fill ($(11)!.4!(21)$) circle (.4mm);
    \fill ($(11)!.5!(21)$) circle (.4mm);
    \fill ($(11)!.6!(21)$) circle (.4mm);
    \fill ($(12)!.4!(22)$) circle (.4mm);
    \fill ($(12)!.5!(22)$) circle (.4mm);
    \fill ($(12)!.6!(22)$) circle (.4mm);

    \draw[-](00)--(01); \draw[-](01)--(02); 
    \draw[-](10)--(11); \draw[-](11)--(12);
    \draw[-](20)--(21); \draw[-](21)--(22); 
    \draw[-](30)--(31); \draw[-](31)--(32); 
  \end{tikzpicture}\ec
  \caption{Coarse structure for a maximal hyperbolic
    isometry.\label{fig:coarse}}
\end{figure}

\begin{defn}[Coarse structure]\label{def:coarse}
  Let $w\in L$ be a hyperbolic euclidean isometry of maximal
  reflection length.  For every $u \in [1,w]^L$ there is a unique $v$
  such that $uv=w$ and we coarsely partition the elements of $[1,w]^L$
  into $3$ rows based on the types of $u$ and $v$ and into columns
  based on the dimensions of their basic invariants.  See
  Figure~\ref{fig:coarse}.  When $u$ or $v$ is hyperbolic, it turns
  out that the other must be an elliptic isometry where every point
  undergoes a motion that is purely horizontal.  In particular, it is
  an elliptic whose fix-set is invariant under vertical translation,
  i.e translation in the direction of the line which is the min-set of
  $w$.  When both $u$ and $v$ are elliptic, it turns out that both
  motions must have non-trivial vertical components and thus neither
  of their fix-sets is invariant under vertical translation.  Within
  each row we grade based on the dimensions of the basic invariants.
  In the bottom row, the dimension of the fix-set of $u$ descreases
  and the dimension of the move-set of $v$ increases as we move from
  left to right.  In the middle row, the dimension of the fix-set of
  $u$ descreases and the dimension of the fix-set of $v$ increases as
  we move from left to right.  And in the top row, the dimension of
  the move-set of $u$ increases and the dimension of the fix-set of
  $v$ increases as we move from left to right.
\end{defn}

The only element in the lower left-hand box is the identity element
corresponding to the factorization $1\cdot w = w$ and the only element
in the upper right-hand box is the element $w$ corresponding to the
factorization $w \cdot 1 = w$.  All other boxes in this grid contain
infinitely many elements.  Descending in the poset order involves
moving to elements in boxes down and/or to the left and covering
relations involve elements in boxes that are adjacent either
vertically or horizontally.  As a consequence, the box an element is
placed in determines its reflection length: its length equals the
number of steps its box is from the lower left-hand corner.

\begin{defn}[Three special boxes]
  There are three particular boxes in this grid that merit additional
  description.  The elements placed in the upper left-hand corner are
  hyperbolic isometries in $[1,w]^L$ of reflection length $2$ which
  means that they are produced by multiplying a pair of reflections
  fixing parallel hyperplanes.  In other words they are pure
  translations $t_\lambda$, and by construction of the model poset
  $P(w)$, the pure translations $T$ which occur in the interval are
  precisely those where the translation vector $\lambda$ is a element
  of the non-linear affine subspace $\mov(w) \subset V$.  The elements
  in the second box in the bottom row have reflection length~$1$,
  i.e. they are themselves reflections and since they are in the
  bottom row, they have fixed hyperplanes invariant under vertical
  translation.  We call this set $R_H$ the \emph{horizontal
    reflections} since they move points in a horizontal direction.
  All such reflections occur in the interval $[1,w]^L$.  And finally
  the first box in the middle row contains reflections whose fixed
  hyperplane is not invariant under vertical translation.  We call
  this set $R_V$ the \emph{vertical reflections} because the motions
  they produce contain a nontrivial vertical component.  All such
  reflections also occur in the interval $[1,w]^L$.
\end{defn}

We conclude our discussion of intervals in the full euclidean isometry
group by describing where in the grid meets and joins fail to exist.

\begin{exmp}[Lattice failure]
  The simplest euclidean isometry whose interval fails to be a lattice
  is a loxodromic ``corkscrew'' motion $w$ in $\R^3$.  This isometry
  has reflection length~$4$ and it has a coarse structure with three
  rows and three columns.  If we consider any pair of hyperbolic
  isometries from the middle box of the top row whose min-sets are
  parallel vertically invariant planes and a pair of horizontal
  reflections in the middle box of the bottom row whose fixed planes
  are parallel to each other and to the min-sets of the chosen
  hyperbolic isometries, then it is straight-forward to check that
  these hyperbolic isometries are distinct minimal upper bounds for
  this pair of elliptic isometries and these elliptic isometries are
  distinct minimal lower bounds for these hyperbolic isometries.  In
  \cite{BrMc-factor} we call this situation a \emph{bowtie}.
\end{exmp}

\part{Crystallographic Garside groups}

In this second part of the article I describe how knowing the
structure of intervals in the full euclidean isometry group leads to
an understanding of similar intervals inside a euclidean Coxeter
group, and how these Coxeter intervals provide the technical
foundation at the heart of our successful attempt to understand
euclidean Artin groups using infinite-type Garside structures.

\section{Coxeter elements}

Let $W = \cox(\wt X_n)$ be an irreducible euclidean Coxeter group
acting geometrically on an $n$-dimensional euclidean space $E$.  The
Coxeter group $W$ is discrete subsgroup of the Lie group $L =
\isom(E)$ and if we continue to view $L$ as a group generated by all
reflections and we view $W$ as the subgroup generated by those
reflections which occur in $W$, then one might naturally expect there
to be a close relationship between the interval $[1,w]^W$ and the
interval $[1,w]^L$ for each $w\in W$.  For generic elements the
connection is not as close as one might hope.  In fact, even the
distance to the origin might be different in the two contexts, which
makes the sets of minimal length reflection factorizations completely
disjoint, as Kyle Petersen and I explored in \cite{McPe10}.  There is,
however, a close connection when $w$ is a Coxeter element of $W$.

\begin{figure}
  \includegraphics{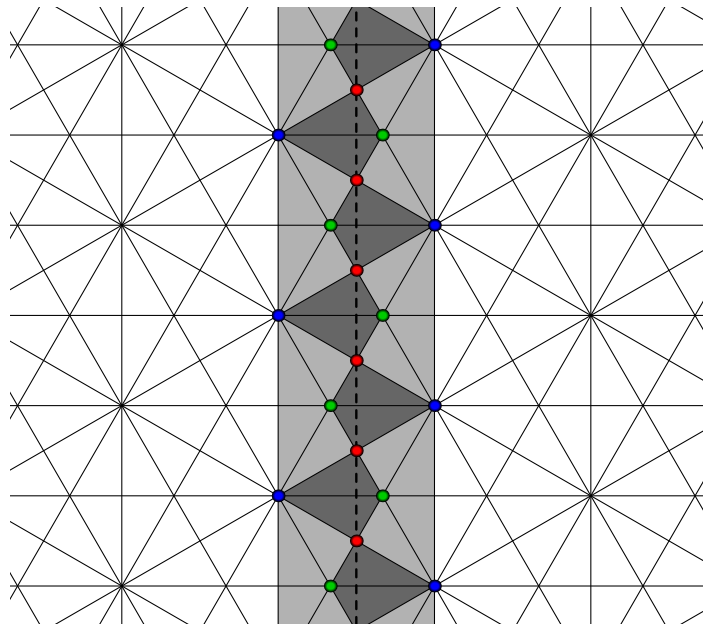}
  \caption{The euclidean Coxeter Group $\cox(\wt G_2)$.\label{fig:g2}}
\end{figure}

\begin{defn}[Coxeter elements]
  Let $W = \cox(\wt X_n)$ be an irreducible euclidean Coxeter group
  with Coxeter generating set $S$.  A \emph{Coxeter element} $w \in W$
  is obtained by multiplying the elements of $S$ in some order.  This
  produces many different Coxeter elements depending on the order in
  which these elements are multiplied, but so long as the diagram $\wt
  X_n$ is a tree, as it is in all cases except for $\wt A_n$, all
  Coxeter elements in $W$ belong to the same conjugacy class and act
  on the corresponding euclidean tiling in the exact same way
  \cite[\McLattGeoEquiv]{Mc-lattice}.  Thus it makes sense to talk
  about \emph{the} Coxeter element in most irreducible euclidean
  contexts.
\end{defn}

Coxeter elements of irreducible euclidean Coxeter groups are
hyperbolic euclidean isometries whose geometric invariants play a
large role in the our understanding of the structure of the interval
$[1,w]^W$.

\begin{defn}[Axial features]
  Let $w$ be a Coxeter element for an irreducible euclidean Coxeter
  group $W=\cox(\wt X_n)$.  It is a hyperbolic isometry whose
  reflection length is $n+1$ when measured in either $W$ or $L$.  In
  $L$ this reflection length is the maximum possible, its min-set
  $\mn(w)$ is a line in $E$ called its \emph{axis} and its move-set
  $\mov(w)$ is a non-linear affine hyperplane in $V$.  The
  top-dimensional simplices whose interior nontrivially intersects the
  axis are called \emph{axial simplices} and the vertices of these
  simplices are \emph{axial vertices}.
\end{defn}

\begin{exmp}[$\cox(\wt G_2)$]
  Figure~\ref{fig:g2} illustrates these ideas for the Coxeter group
  $\cox(\wt G_2)$.  Its Coxeter element is a glide reflection whose
  glide axis, i.e. its min-set, is shown as a dashed line.  The
  corresponding axial simplices are heavily shaded their axial
  vertices are shown as enlarged dots.
\end{exmp}

For Coxeter elements, the interval $[1,w]^W$ is a restriction of
edge-labeled subposet $[1,w]^L$ to the union of minimal length paths
in the Cayley graph of $L$ from $v_1$ to $v_w$ where every edge is
labeled by a reflection in $W = \cox(\wt X_n)$.  The original product
of elements in the Coxeter generating set $S$ which produces $w$ as a
Coxeter element is one such minimal length path.  As such, the
elements of this \emph{Coxeter interval} $[1,w]^W$ have a coarse
structure as described in Definition~\ref{def:coarse}.  The first
difference we find is that whereas every reflection in $L$ labels some
edge in the interval $[1,w]^L$, in the Coxeter interval $[1,w]^W$ only
a proper subset of the reflections in $W$ actually label edges in the
interval.  For the unused reflections do not occur in a minimal length
factorization of $w$ where every reflection must belong to $W$.  In
\cite[\McLattRefl]{Mc-lattice} the reflections that do occur as edge
labels in the interval are precisely characterized as follows.

\begin{thm}[Reflection generators]
  Let $w$ be a Coxeter element of an irreducible euclidean Coxeter
  group $W=\cox(\wt X_n)$. A reflection labels an edge in the
  interval $[1,w]^W$ iff its fixed hyperplane contains an axial
  vertex.
\end{thm}

If we separate the reflections labeling edges in $[1,w]^W$ into those
which are horizontal and those which are vertical, in the sense
defined in the previous section, then there are infinitely many
vertical reflections (all those whose hyperplanes cross the Coxeter
axis) and a finite number of horizontal reflections (those whose
hyperplanes bound the convext hull of the axial simplices).  More
generally, the coarse structure of the interval $[1,w]^W$ will have
only finitely many elements in each box along the top row and finitely
many elements in each box along the bottom row.  The boxes in the
middle row, on the other hand, have infinitely many elements in each.
Nevertheless, there is a periodicity to the convex hull of the axial
simplices and this means that the infinitely many elements in each box
in the middle row falls into a finite number of infinitely repeating
patterns.  We illustrate this with the $\wt G_2$ Coxeter group where
one can view the entire euclidean tiling.

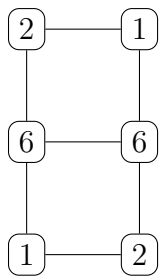
\begin{figure}
  \bc\begin{tikzpicture}[node distance=15mm, auto]
  \tikzstyle{every node}=[rounded corners,draw] \node(00){$1$};
  \node[right of=00](10){$2$};

    \node[above of=00](01){$6$}; 
    \node[right of=01](11){$6$}; 

    \node[above of=01](02){$2$}; 
    \node[right of=02](12){$1$}; 

    \draw[-](00)--(10); 
    \draw[-](01)--(11); 
    \draw[-](02)--(12); 

    \draw[-](00)--(01); \draw[-](01)--(02); 
    \draw[-](10)--(11); \draw[-](11)--(12);
  \end{tikzpicture}\ec
  \caption{Coarse structure for the $\wt G_2$
    interval.\label{fig:g2-coarse}}
\end{figure}

\begin{exmp}[Coarse structure of the $\wt G_2$ interval]
  As can be seen in Figure~\ref{fig:g2}, the euclidean Coxeter group
  $W=\cox(\wt G_2)$ has exactly $2$ horizontal reflections, the ones
  with vertical fixed lines which bound the lightly shaded region and
  this is indicated by the $2$ in the second box in the bottom row of
  the coarse structure for the interval $[1,w]^W$ shown schematically
  in Figure~\ref{fig:g2-coarse}. On the other hand, there are $6$
  essentially different ways that a fixed line can cross the glide
  axis and the corresponding $6$ infinite families of vertical
  reflections below $w$ are indicated by the $6$ in the first box of
  the middle row.  Similarly, there are $6$ infinite families of
  rotations fixing an axial vertex represented by the $6$ in the
  second box of the middle row and exactly two pure translations below
  $w$ indicated by the $2$ in the first box of the top row.  Finally,
  although it is not immmediately obvious, it turns out that this
  bounded, graded, self-dual poset with finite height and an infinite
  number of elements is a lattice.
\end{exmp}

For a more interesting example, consider the largest sporadic
irreducible euclidean Coxeter group $\cox(\wt E_8)$.

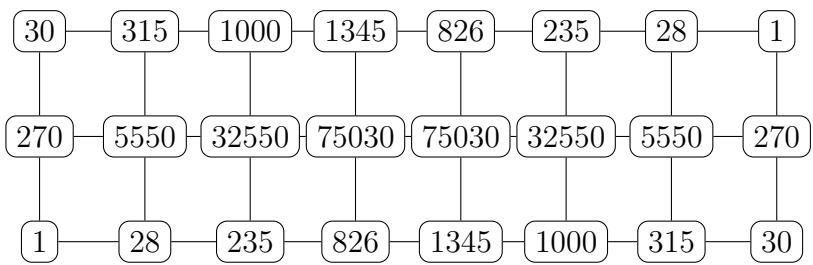
\begin{figure}
  \bc\begin{tikzpicture}[node distance=14mm, auto]
    \tikzstyle{every node}=[rounded corners,draw]
    \node(00){$1$}; 
    \node[right of=00](10){$28$}; 
    \node[right of=10](20){$235$}; 
    \node[right of=20](30){$826$}; 
    \node[right of=30](40){$1345$}; 
    \node[right of=40](50){$1000$}; 
    \node[right of=50](60){$315$}; 
    \node[right of=60](70){$30$}; 

    \node[above of=00](01){$270$}; 
    \node[right of=01](11){$5550$}; 
    \node[right of=11](21){$32550$}; 
    \node[right of=21](31){$75030$}; 
    \node[right of=31](41){$75030$}; 
    \node[right of=41](51){$32550$}; 
    \node[right of=51](61){$5550$}; 
    \node[right of=61](71){$270$}; 

    \node[above of=01](02){$30$}; 
    \node[right of=02](12){$315$}; 
    \node[right of=12](22){$1000$}; 
    \node[right of=22](32){$1345$}; 
    \node[right of=32](42){$826$}; 
    \node[right of=42](52){$235$}; 
    \node[right of=52](62){$28$}; 
    \node[right of=62](72){$1$}; 

    \draw[-](00)--(10); \draw[-](10)--(20); \draw[-](20)--(30); \draw[-](30)--(40); 
    \draw[-](40)--(50); \draw[-](50)--(60); \draw[-](60)--(70); 

    \draw[-](01)--(11); \draw[-](11)--(21); \draw[-](21)--(31); \draw[-](31)--(41); 
    \draw[-](41)--(51); \draw[-](51)--(61); \draw[-](61)--(71); 

    \draw[-](02)--(12); \draw[-](12)--(22); \draw[-](22)--(32); \draw[-](32)--(42); 
    \draw[-](42)--(52); \draw[-](52)--(62); \draw[-](62)--(72); 

    \draw[-](00)--(01); \draw[-](01)--(02); 
    \draw[-](10)--(11); \draw[-](11)--(12);
    \draw[-](20)--(21); \draw[-](21)--(22);
    \draw[-](30)--(31); \draw[-](31)--(32);
    \draw[-](40)--(41); \draw[-](41)--(42);
    \draw[-](50)--(51); \draw[-](51)--(52);
    \draw[-](60)--(61); \draw[-](61)--(62);
    \draw[-](70)--(71); \draw[-](71)--(72);
  \end{tikzpicture}\ec
  \caption{Coarse structure for the $\wt E_8$ interval.\label{fig:e8-coarse}}
\end{figure}

\begin{exmp}[Coarse structure of the $\wt E_8$ interval]\label{ex:e8}
  Let $w$ be a Coxeter element for the Coxeter group $W = \cox(\wt
  E_8)$.  The coarse structure of the interval $[1,w]^W$ is shown
  schematically in Figure~\ref{fig:e8-coarse}.  As in the case of
  $\cox(\wt G_2)$, the numbers listed in the top row and in bottom row
  indicate the actual number of elements in each box, but the numbers
  in the middle row only indicate the number of infinite families of
  such elements.  The equivalence relation used is that two middle row
  elements below $w$ belong to the same family iff they differ by
  (conjugation by) a translation of the tiling in the direction of the
  Coxeter axis.  Thus, from the coarse structure we see that the
  interval $[1,w]^W$ has $28$ horizontal reflections (second box in
  the bottom row), $30$ translations (first box in the top row), $270$
  infinite families of vertical reflections (first box in the middle
  row) and $5550$ infinite families of elements that rotation around a
  $6$-dimensional fix-set that is not invariant under vertical
  translation (second box in the middle row), and so on.
  Representatives of each family were computed using a program I wrote
  called \texttt{euclid.sage} that is available from my webpage.  The
  software also checks whether this bounded, graded, self-dual poset
  of finite height with an infinite number of elements is a lattice
  and in this case the answer is ``No''.
\end{exmp}

\section{Dual Artin groups\label{sec:dual-artin}}

Now that the technical foundations are in place, it is time to shift
our attention to the irreducible euclidean Artin groups themselves, to
finally explain why intervals in euclidean Coxeter groups are relevant
and why we are interested in whether or not these intervals are
lattices.  The answers are relatively straightforward.  First,
intervals in irreducible euclidean Coxeter groups can be used to give
alternative, so called ``dual'' presentations for irreducible
euclidean Artin groups as I show in \cite{Mc-lattice}.  Next, when
these Coxeter intervals are lattices, the dual Artin group has an
infinite-type Garside structure and groups with Garside structures
have good computational properties.  This ``grand scheme'', closely
related to the approach taken by Fran\c{c}ois Digne in \cite{Digne06}
and \cite{Digne12} was the initial strategy by which my coauthors and
we hoped to use it to understand arbitary euclidean Artin groups.
Unfortunately, a detailed examination of the groups themselves caused
this scheme to fail because the Coxeter intervals turned out to be
more poorly behaved than expected.  Nonetheless, a modified grand
scheme, described in the later sections, does eventually succeed.  The
first step is to understand how intervals lead to presentations of new
groups.

\begin{defn}[Interval groups]
  Let $G$ be a group with a fixed symmetric discretely weighted
  generating set and let $[1,g]^G$ be an interval in $G$, viewed as an
  edge-labeled directed graph sitting inside the Cayley graph of $G$.
  The \emph{interval group} $G_g$ is a new group generated by the
  labels of edges in the interval subject to the set of all relations
  that are visible in the interval.  Since there is a natural function
  from the generators of $G_g$ to $G$ and since the set of relations
  used to define $G_g$ is a subset of the relations which hold in $G$,
  this function extends to a group homomorphism $G_g \to G$.  If,
  moreover, the labels on the edges in the interval $[1,g]^G$ include
  a generating set for $G$ then this natural map is onto.
\end{defn}

To see how this works in practice, consider the following example.

\begin{exmp}[Noncrossing partitions]
  Let $G = \sym_n$ be the symmetric group on $n$ elements, let $g$ be
  the $n$-cycle $(1,2,\ldots, n)$ and fix the full set of
  tranpositions as its generating set.  It turns out in this case that
  the poset structure of the interval $[1,g]^G$ is a well-known
  combinatorial object called the \emph{noncrossing partition lattice}
  defined as follows.  Start with a convex regular $n$-gon in the
  plane whose vertices are labeled $1$ through $n$ in a clockwise
  fashion.  A partitioning of its vertex set $[n]=\{1,2,\ldots,n\}$ is
  called \emph{noncrossing} if the convex hulls of the blocks of the
  partition are pairwise disjoint.  For example, the partition $\{
  \{1,3\},\{2,4\}\}$ of $[4]$ is not a noncrossing partition because
  the convex hulls are two line segments that intersect.  Noncrossing
  partitions can be ordered by declaring $\sigma < \tau$ when every
  block of $\sigma$ is a subset of a block of $\tau$ and a noncrossing
  partition can be converted into a permutation by clockwise permuting
  the vertices in the boundary of the convex hull of each block.  This
  function defines a poset isomorphism between the noncrossing
  partition lattice $NC_n$ and the Coxeter interval $[1,g]^G$.  See
  \cite{Mc06} for an elementary discussion of these ideas.  When
  $n=3$, there are only $3$ transpositions, the noncrossing parititon
  lattice is particularly simple and the presentation for the
  corresponding interval group $G_g$ is $\langle a,b,c \mid ab=bc=ca
  \rangle$, which is an alternate presentation for the $3$-string
  braid group.
\end{exmp}

This example, and its generalization to all spherical Coxeter groups
described below, leads to the following general definition.

\begin{defn}[Dual Artin groups]
  Let $W=\cox(\Gamma)$ be an arbitrary Coxeter group viewed as a group
  generated by its full set of reflections and let $w$ be one of its
  Coxter elements.  The interval group $W_w$ defined by the interval
  $[1,w]^W$ is called a \emph{dual Artin group} and denoted
  $\dart(\Gamma,w)$.  The notation is meant to highlight the fact that
  in general Coxeter groups there are geometrically distinct Coxeter
  elements and thus there are distinct dual presentations which
  heavily depend on the choice of Coxeter element.
\end{defn}

The study of dual presentations in general is motivated by the work of
Davis Besis \cite{Be03} and Tom Brady and Colum Watt \cite{BrWa02a} on
spherical Artin groups.  Here are their main results translated into
this terminology.

\begin{thm}[Dual spherical Artin groups]
  If $W =\cox(X_n)$ is a spherical Coxeter group generated by its
  reflections, and $w$ is a Coxeter element of $W$, then the Coxeter
  interval $[1,w]^W$ is isomorphic to the $W$-noncrossing partition
  lattice and the interval group $W_w =\dart(X_n,w)$ is isomorphic to
  the corresponding Artin group $\art(X_n)$.
\end{thm}

The terminology ``dual Artin group'' was introduced in
\cite{Mc-lattice} because in general it is not known whether or not
Artin groups and dual Artin groups are isomorphic.  Fortunately, I was
able to establish that they are isomorphic in the euclidean case
\cite{Mc-lattice}.

\begin{thm}[Dual euclidean Artin groups]
  If $W =\cox(\wt X_n)$ is an irreducible euclidean Coxeter group
  generated by its reflections, and $w$ is a Coxeter element, then the
  dual Artin group $W_w = \dart(\wt X_n,w)$ is naturally isomorphic to
  $\art(\wt X_n)$.
\end{thm}

The proof uses a result from quiver representation theory to greatly
simplify the dual presentations for dual euclidean Artin groups and it
is then relatively straight-forward to establish that these simplified
presentations define the same groups as do the standard euclidean
Artin presentations \cite{Mc-lattice}.  In other words, the interval
$[1,w]^W$ gives a new presentation (with infinitely many generators
and infinitely many relations) of the corresponding Artin group in the
irreducible euclidean case.  The interest in whether or not these
Coxeter intervals are lattices has to do with the following result
which is essentially due to David Bessis \cite{Be03}.  See
\cite{Mc-lattice} for a more detailed discussion.

\begin{thm}[Sufficient conditions]
  Let $G$ be a group with a fixed symmetric discretely weighted
  generating set closed under conjugation.  For each $g\in G$, if the
  interval $[1,g]^G$ is a lattice, then $G_g$ is a Garside group,
  possibly of infinite-type.
\end{thm}

In this article, Garside structures are treated as a black box.  For a
book-length discussion of Garside structures see
\cite{DDGKM-garside}. The main idea is that if there is a portion of
the Cayley graph which contains a generating set and has well-defined
meets and joins (plus a few more technical conditions) then this local
lattice structure can be used to systematically construct normal forms
for all group elements, thereby solving the word problem for the group
and allowing one to construct a finite dimensional classifying space.

\begin{thm}[Garside consequences]
  If an interval group $G_g$ is a Garside group, possibly of
  infinite-type, then $G_g$ is a torsion-free group with a solvable
  word problem and a finite dimensional classifying space.
\end{thm}

\begin{rem}[Artin groups with Garside structures]
  By the early years of the new millenium, many dual Artin groups were
  known to have Garside structures.  We have already mentioned that
  the dual Artin group $\art(\Gamma)$ is Garside when $\cox(\Gamma)$
  is spherical (due to Bessis \cite{Be03} and Brady-Watt
  \cite{BrWa02a} independently) and when $\Gamma$ is an extended
  Dynkin diagram of type $\wt A_n$ or $\wt C_n$ (due to Digne
  \cite{Digne06,Digne12}).  David Bessis also proved this for the free
  group, thought of as the Artin group where every $m_{ij}=\infty$
  \cite{Be06}.  In addition, there are unpublished results due to
  myself and John Crisp which show this to be true for all Artin
  groups with at most $3$ standard generators and also for all Artin
  groups defined by a diagram in which every $m_{ij}$ is at least $6$.
  Note that the $3$ standard generator result means that the three
  planar Artin groups investigated by Craig Squier in \cite{Squier87}
  have dual presentations which are Garside structures (as do all of
  the $3$-generator Artin groups whose Coxeter groups naturally act on
  the hyperbolic plane).  In particular, $\art(\wt G_2)$ has an
  infinite-type Garside structure.
\end{rem}

This list of results naturally lead one to conjecture that Coxeter
intervals using Coxeter elements are always lattices and that all the
corresponding dual Artin groups are Garside groups.  Unfortunately, as
we have already seen in Example~\ref{ex:e8}, this natural conjecture
turns out to be too optimistic and false even for some of the
irreducible euclidean Coxeter groups such as the group $W = \cox(\wt
E_8)$.

\section{Horizontal roots}

Several years ago John Crisp and I systematically investigated whether
every irreducible euclidean Coxeter group has a Coxeter interval which
is a lattice, and what we found was not what we expected to find.  The
only irreducible euclidean Coxeter groups whose Coxeter intervals are
lattices are those of type $\wt A_n$, $\wt C_n$, and $\wt G_2$.  In
other words, the only dual Artin groups with Garside structures are
those where the group structure was already well understood by Craig
Squier and/or Fran\c{c}ois Digne.  Further investigation revealed that
the reason why a euclidean Coxeter group might have a Coxeter interval
which failed to be a lattice is closely related to the structure of
what we call its horizontal root system.

\begin{defn}[Horizontal root system]
  Let $W = \cox(\wt X_n)$ be an irreducible euclidean Coxeter group
  and let $w$ be a Coxeter element of $W$.  The horizontal reflections
  in the interval $[1,w]^W$ are those whose roots are orthogonal to
  the direction of the Coxeter axis of $w$.  The set of all such roots
  form a subroot system of the full root system of $W$ that we call
  the \emph{horizontal root system}.
\end{defn}  

\begin{table}
  \begin{center}
    $\begin{array}{|c|l|}
      \hline
      \textrm{Type} & \textrm{Horizontal root system}\\
      \hline
      A_n & \Phi_{A_{p-1}} \cup \Phi_{A_{q-1}} \\
      C_n & \Phi_{A_{n-1}} \\
      B_n & \Phi_{A_1} \cup \Phi_{A_{n-2}} \\
      D_n & \Phi_{A_1} \cup \Phi_{A_1} \cup \Phi_{A_{n-3}} \\
      \hline
      G_2 & \Phi_{A_1} \\
      F_4 & \Phi_{A_1} \cup \Phi_{A_2} \\
      E_6 & \Phi_{A_1} \cup \Phi_{A_2} \cup \Phi_{A_2} \\
      E_7 & \Phi_{A_1} \cup \Phi_{A_2} \cup \Phi_{A_3} \\
      E_8 & \Phi_{A_1} \cup \Phi_{A_2} \cup \Phi_{A_4} \\
      \hline
    \end{array}$
  \end{center}\vspace*{1em}
  \caption{Horizontal root systems by type.\label{tbl:horizontal}}
\end{table}

It turns out that horizontal root system is easy to describe as a
subdiagram of the original extended Dynkin diagram.

\begin{rem}[Finding horizontal roots]
  Let $W = \cox(\wt X_n)$ be an irreducible euclidean Coxeter group
  with Coxeter element $w$.  The horizontal root system with respect
  to $w$ is itself a root system for a spherical Coxeter group whose
  Dynkin diagram is obtained by removing two dots from the extended
  Dynkin diagram $\wt X_n$ or one dot from Dynkin diagram $X_n$.  In
  Figure~\ref{fig:dynkin} the dots to be removed are slightly larger
  than the others.  Removing the large white dot produces the Dynkin
  diagram of type $X_n$.  Also removing the large shaded dot produces
  the diagram for the horizontal root system.  In all cases, the
  shaded dot is the long end of a multiple bond or the branch point if
  either exists in $X_n$.  The only case where neither exists is in
  type $\wt A_n$, where the shaded dot might be any of the remaining
  dots and different choices arise from the geometrically different
  choices for the Coxeter element in this case. 
\end{rem}

The types of the horizontal root systems are listed in
Table~\ref{tbl:horizontal}.  The key property turns out ot be whether
or not the remaining diagram is connected, or equivalently, whether or
not the horizontal root system is reducible.  The horizontal root
systems for types $C$ and $G$ are irreducible, the horizontal root
systems for types $B$, $D$, $E$ and $F$ are reducible and for type $A$
it depends on the choice of Coxeter element.  The following theorem is
a restatement of and explanation for the computational result
originally discovered in collaboration with John Crisp.  It is proved
in \cite{Mc-lattice}, an article which is, morally speaking, the
result of a collaboration with John Crisp even though it was only
written up after he left research mathematics.
  
\begin{thm}[Failure of the lattice property]
  The interval $[1,w]^W$ is a lattice iff the horizontal root system
  is irreducible.  In particular, types $C$ and $G$ are lattices,
  types $B$, $D$, $E$ and $F$ are not, and for type $A$ the answer
  depends on the choice of Coxeter element.
\end{thm}

As a consequence of this theorem, it is clear which irreducible
euclidean Artin groups have dual presentations that lead to Garside
structures.

\begin{cor}
  The dual euclidean Artin group $\dart(\wt X_n,w)$ is Garside when
  $X$ is $C$ or $G$ and it is not Garside when $X$ is $B$, $D$, $E$ or
  $F$.  When the group has type $A$ there are distinct dual
  presentations and the one investigated by Digne is the only one that
  is Garside.
\end{cor}

At this point, the grand scheme has failed and no additional euclidean
Artin groups have been understood.  There is, however, a ray of hope.
The explicit nature of the euclidean model posets enables an explicit
examination of the pairs of elements which fail to have well-defined
meets or well-defined joins.  It turns out that pairs of elements with
no well-defined join must occur in the bottom row of the coarse
structure and pairs of elements with no well-defined meet must occur
in the top row of the coarse structure.  Since the top and bottom rows
only contain finitely many elements, this means that out of the
infinitely many pairs of elements in the infinite interval $[1,w]^W$,
only finitely pairs fail to be well-behaved.  This leaves open the
possibility that one can systematically fix these finitely many
failures.

\section{New groups}

Let $W =\cox(\wt X_n)$ be an irreducible euclidean Coxeter group with
Coxeter element $w$.  As remarked above, the finitely pairs of
elements in the Coxeter interval $[1,w]^W$ which cause the lattice
property to fail all occur in the top and bottom rows of its coarse
structure.  Thus it makes sense to focus on the subgroup corresponding
to this portion of the interval.  Its structure is closely related to
an elementary group which does not appear to have a standard name in
the literature.  In fact, we have not found any references to it in
the literature so far.

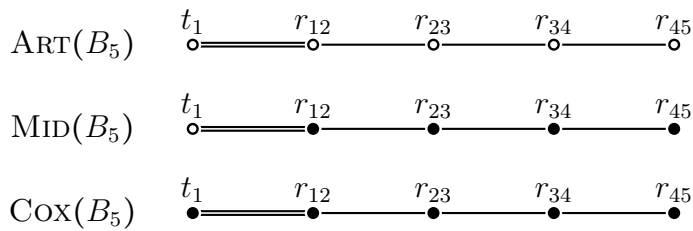
\begin{figure}
  \bc\begin{tikzpicture}[scale=1.6]
    \foreach \y in {-.7,0,.7} {
      \draw[double,thick] (-2,\y)--(-1,\y);
      \draw[-,thick] (-1,\y)--(2,\y);
      \foreach \x in {-2,-1,0,1,2} {
        \fill[color=white] (\x,\y) circle (.7mm);
        \fill (\x,\y) circle (.5mm);}
      \draw (-2,\y) node [anchor=south] {$t_1$};
      \draw (-1,\y) node [anchor=south] {$r_{12}$};
      \draw (0,\y) node [anchor=south] {$r_{23}$};
      \draw (1,\y) node [anchor=south] {$r_{34}$};
      \draw (2,\y) node [anchor=south] {$r_{45}$};
    }
    \node at (-3,.7) {$\art(B_5)$};
    \node at (-3,0) {$\midd(B_5)$};
    \node at (-3,-.7) {$\cox(B_5)$};
    \foreach \x in {-2,-1,0,1,2} \fill[color=white] (\x,.7) circle (.3mm);
    \fill[color=white] (-2,0) circle (.3mm);
  \end{tikzpicture}\ec
  \caption{Presentation diagrams for the Artin, middle and Coxeter
    groups of type $B_5$.\label{fig:mid-dia}}
\end{figure}

\begin{defn}[Middle groups]
  Consider the group of isometries of $\Z^n$ in $\R^n$ generated by
  coordinate permutations and integral translations.  We call this
  group the \emph{middle group} and denote it $\midd(B_n)$.  It is
  generated by the reflections $r_{ij}$ that switch coordinates $i$
  and $j$ and the translations $t_i$ that adds $1$ to the $i$-th
  coordinate and it is a semidirect product $\Z^n \rtimes \sym_n$ with
  the translations $t_i$ generating the normal free abelian subgroup.
  A standard minimal generating set for $\midd(B_n)$ is the set
  $\{t_1\} \cup \{r_{12}, r_{23}, \ldots, r_{n-1 n}\}$ and it has a
  presentation similar to $\art(B_n)$ and $\cox(B_n)$
  \cite{McSu-artin-euclid}.  See Figure~\ref{fig:mid-dia}.  A solid
  dot means the corresponding generator has order $2$ and an empty dot
  means the corresponding generator has infinite order.
\end{defn}

If we consider $\midd(B_n)$ as a group generated by the full set of
translations and reflections, then the factorizations of $w=t_1 r_{12}
r_{23} \cdots r_{n-1 n}$ form an interval isomorphic to the type~$B$
noncrossing partition lattice, exactly the same poset as the Coxeter
interval in the spherical Coxeter group $\cox(B_n)$.  This explains
the use of $B_n$ in the notation.  The name ``middle group'' is
suggested by its location in the center of a diagram that shows its
relation to several closely related Coxeter groups and Artin groups.
See Figure~\ref{fig:mid-rel}.  The top row is the short exact sequence
that is often used to understand $\art(\wt A_{n-1})$.  Geometrically
middle groups are easy to recognize as a symmetric group generated by
reflections and a translation with a component out of this subspace.
Middle groups are introduced in \cite{McSu-artin-euclid} in order to
succinctly describe the structure of the diagonal subgroup build from
the top and bottom rows of the coarse structure of a euclidean Coxeter
interval.

\begin{figure}
  \begin{center}
    $\begin{array}{ccccc}
      \art(\wt A_{n-1}) & \into & \art(B_n) & \onto & \Z\\
      \twoheaddownarrow & & \twoheaddownarrow & & \twoheaddownarrow\\
      \cox(\wt A_{n-1}) & \into & \midd(B_n) & \onto & \Z\\
      & & \twoheaddownarrow\\
      & & \cox(B_n)\\
    \end{array}$
  \end{center}
  \caption{Relatives of middle groups.\label{fig:mid-rel}}
\end{figure}

\begin{defn}[Diagonal subgroup]
  Let $W = \cox(\wt X_n)$ be an irreducible euclidean Coxeter group
  with Coxeter element $w$ and let $R_H$ and $T$ denote the set of
  horizontal reflections and translations contained in the Coxeter
  interval $[1,w]^W$.  In addition let $D$ be the subgroup of $W$
  generated by the set $R_H \cup T$.  If we assign a weight of $1$ to
  each horizontal reflection and a weight of $2$ to each translation,
  then distances in the Cayley graph of $D$ match distances in the
  Cayley graph of $W$ and the interval $[1,w]^D$ is an induced
  subposet of the Coxeter interval $[1,w]^W$ consisting of only the
  top and bottom rows.  We write $D_w$ to denote the interval group
  defined by this restricted interval.
\end{defn}

The interal $[1,w]^D$ is almost a direct product of posets and the
group $D$ is almost a direct product of middle groups.  More
precisely, the poset $[1,w]^D$ is almost a direct product of type~$B$
noncrossing partitions lattices and the missing elements are added if
we factor the translations of $D$.

\begin{defn}[Factored translations]
  Each pure translation $t$ in $[1,w]^D$ projects nontrivially onto
  the Coxeter axis and onto each of the irreducible components of the
  horizontal root system of the corresponding Coxeter group $W =
  \cox(\wt X_n)$.  Let $t^{(i)}$ be the translation which agrees with
  $t$ on the $i$-th component and contains $1/k$ of the translation in
  the Coxeter direction where $k$ is the number of irreducible
  components of the horizontal root system.  We call each $t^{(i)}$ a
  \emph{factored translation} and let $T_F$ denote the set of all such
  factored translations.
\end{defn}

We use the factored translations to introduce a slightly larger group.

\begin{defn}[Factorable groups]
  The \emph{factorable group} $F$ is defined as the group of euclidean
  isometries generated by $R_H \cup T_F$.  It is crystallographic in
  the sense that it acts geometrically on euclidean space but it is
  not a Coxeter group in general because it is not generated by
  reflections.  If we assign a weight of $2/k$ to each factored
  translation then distances in the Cayley graph of $F$ agree with
  those in $D$ and the interval $[1,w]^D$ is an induced subinterval of
  $[1,w]^F$.  The main advantage of the interval $[1,w]^F$ is that it
  factors as a direct product of type~$B$ noncrossing partition
  lattices with one factor for each irreducible component of the
  horzontal root system.  The edge labels in the $i$-th factor poset
  correspond to the factored translations which project nontrivially
  onto $i$-th component of the horizontal root system together with
  the horizontal reflections whose roots belong to this component.
  Moreover, the isometries that occur as labels in any one factor
  generator a group isomorphic to a middle group acting on the
  subspace whose directions are spanned by a component of the
  horizontal root system and the direction of the Coxeter axis.  The
  structure of the factorable group $F$ is not quite as a clean since
  each of the component middle groups contain central elements which
  are pure translations in the direction of the Coxeter axis.  Thus
  $F$ is merely a central product of the associated middle groups.
\end{defn}

In addtion to the groups $W$, $D$, and $F$, we introduce several other
groups that help to clarify the structure of the corresponding
euclidean Artin group.

\begin{table}
  \begin{tabular}{|c|c|l|}
    \hline
    \textbf{Name} & \textbf{Symbol} & \textbf{Generating set}\\
    \hline
    Horizontal & $\chor$ & $R_H$\\
    Diagonal & $\cdiag$ & $R_H \cup T$\\
    Coxeter & $\ccox$ & $R_H \cup R_V\ (\cup\ T )$\\
    Factorable & $\cfac$ & $R_H \cup T_F\ (\cup\ T)$\\
    Crystallographic & $\ccryst\ $ & $R_H \cup R_V \cup T_F\ (\cup\ T)$\\
    \hline
  \end{tabular}\vspace*{1em}
  \caption{Five generating sets.\label{tbl:gen-sets}}
\end{table}

\begin{figure}
  \bc\begin{tikzpicture}[node distance=2cm, auto]
  \node(GH){$\ghor$}; 
  \node[right of=GH](GD){$\gdiag$}; 
  \node[right of=GD](GA){$\gart$}; 
  \node(GF)[right of=GD, above of=GD, node distance=1cm]{$\gfac$}; 
  \node[right of=GF](G){$\ggar$};
  \node[below of=GF](CF){$\cfac$}; 
  \node[below of=G](C){$\ccryst$}; 
  \node[below of=GH](CH){$\chor$};
  \node[below of=GD](CD){$\cdiag$}; 
  \node[below of=GA](CC){$\ccox$};
  \draw[->](GF)--(G); 
  \draw[->](GH)--(GD); 
  \draw[->](GD)--(GA); 
  \draw[->](GD)--(GF); 
  \draw[->](GA)--(G);
  \draw[right hook->](CF)--(C); 
  \draw[right hook->](CH)--(CD); 
  \draw[right hook->](CD)--(CC);
  \draw[right hook->](CD)--(CF); 
  \draw[right hook->](CC)--(C);
  \draw[->>](GF)--(CF); 
  \draw[->>](G)--(C); 
  \draw[->>](GH)--(CH);
  \draw[->>](GD)--(CD); 
  \draw[->>](GA)--(CC);
  \end{tikzpicture}\ec
  \caption{Ten groups.\label{fig:ten-gps}}
\end{figure}
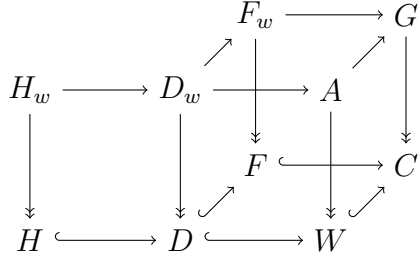

\begin{defn}[Ten groups]
  Let $W=\cox(\wt X_n)$ be an irreducible euclidean Coxeter group with
  Coxeter element $w$ and let $D$ and $F$ be the diagonal and
  factorable groups acting on $n$-dimensional euclidean space as
  defined above.  There are two other euclidean isometry groups we
  need to define.  Let $H$ denote the subgroup of $W$ generated by
  horizontal reflections $R_H$ that label edges in the interval
  $[1,w]^W$ and let $C$ be the group generated by the set $R_H \cup
  R_V \cup T_F (\cup T)$ of all generating isometries considered so
  far.  The set $T$ of translations is in parentheses because its
  elements are products of other generators.  Note that $C$, like $F$,
  is crystallographic in that it acts geometrically on euclidean space
  but it is not in general a Coxeter group since we have added
  translation generators which are not products of reflections
  contained in $W$.  The five generating sets introduced are listed in
  Table~\ref{tbl:gen-sets} along with the euclidean isometry groups
  they generate.  Of these five groups $H$ and $W$ are Coxeter groups,
  while $D$, $F$ and $C$ are merely crystallographic.
  
  Next we construct five groups defined by presentations.  Let $D_w$,
  $F_w$, $W_w$ and $C_w$ denote the interval groups defined by the
  interval $[1,w]$ in each of the four contexts, but note that we
  write $A = W_w$ and $G = C_w$ since these turn out to be the
  corresponding Artin group and a previously unstudied Garside group,
  respectively.  A final group $H_w$ is defined by a presentation with
  $R_H$ as its generators and subject to the relations among these
  generators which are visible in the interval $[1,w]^W$.  This is not
  quite an interval group since there is no interval of the form
  $[1,w]^H$.  This is because the element $w$ is not in the subgroup
  $H$ as it requires a vertical motion in order to be constructed.
  Some of the maps between these ten groups are shown in
  Figure~\ref{fig:ten-gps}.  The maps in the bottom level are the
  natural inclusions among these five euclidean isometry groups, the
  vertical arrows are the projections from the five groups defined by
  presentations to the groups from which they were constructed, and
  the maps in the top level are the natural homomorphisms that extend
  the inclusions on their generating sets.
\end{defn}

A review of some of the various posets and groups associated with the
euclidean Coxeter group of type $\wt E_8$ might help to clarify these
definitions.

\begin{figure}
  \bc\begin{tikzpicture}
    \begin{scope}
      \node at (-1,0) {$\wt E_8$};
      \drawEdge{(0,0)}{(6,0)}
      \drawDashEdge{(6,0)}{(6,1)}
      \drawEdge{(2,0)}{(2,1)}
      \foreach \x in {0,1,3,4,5,6} {\drawRegDot{(\x,0)}}
      \drawRegDot{(2,1)}
      \drawESpDot{(2,0)}
      \drawSpeDot{(6,1)}
    \end{scope}
  \end{tikzpicture}\ec
  \caption{The $\wt E_8$ diagram.\label{fig:e8-dia}}
\end{figure}
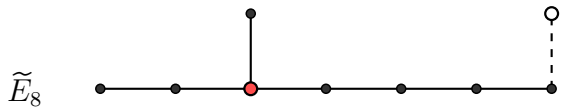

\begin{exmp}[Groups of type $\wt E_8$]
  When $W$ is the irreducible euclidean Coxeter group of type $\wt
  E_8$, its horizontal root system decomposes as $\Phi_{A_1} \cup
  \Phi_{A_2} \cup \Phi_{A_4}$.  See Figure~\ref{fig:e8-dia}.  The
  factorable group $F$ is a central product of middle groups
  $\midd(B_2)$, $\midd(B_3)$ and $\midd(B_5)$.  The interval $[1,w]^F$
  is isomorphic to the direct product $NC_{B_2} \times NC_{B_3} \times
  NC_{B_5}$ of type $B$ noncrossing partition lattices, and the
  interval $[1,w]^F$ defines an interval group $\gfac$ which is a
  direct product $\art(B_2) \times \art(B_3) \times \art(B_5)$ of
  type~$B$ spherical Artin groups.  The horizontal reflections in the
  interval $[1,w]^W$ generate a group $\chor$ isomorphic to a direct
  product $\cox(\wt A_1) \times \cox(\wt A_2) \times \cox(\wt A_4)$ of
  type~$A$ euclidean Coxeter groups, and the relations among these
  reflections visible in the interval $[1,w]^W$ define a group $\ghor$
  which is isomorphic to a direct product $\art(\wt A_1) \times
  \art(\wt A_2) \times \art(\wt A_4)$ of type~$A$ euclidean Artin
  groups.
\end{exmp}

\section{Structural results}

The addition of the factored translations as extra generators
completely solves the lattice problem and makes it possible to to
prove the three main structural results established in by myself and
Rob Sulway in \cite{McSu-artin-euclid}.  The first establishes
the existence of a new class of Garside groups based on intervals in
the crystallographic groups $C$ introduced in the previous section.

\setcounter{main}{0}
\begin{main}[Crystallographic Garside groups]
  If $C = \cryst(\wt X_n,w)$ is the crystallographic group obtained by
  adding the factored translations to the generating set of the
  irreducible euclidean Coxeter group $W = \cox(\wt X_n)$, then the
  interval $[1,w]^C$ in the Cayley graph of $C$ is a lattice.  As a
  consequence, this interval defines an group $G = C_w$ with an
  infinite-type Garside structure.
\end{main}

As is typical, the most difficult step in the entire article is the
proof that these augmented intervals are lattices.  Using the program
\texttt{euclid.sage}, we verify that this is the case for all
irreducible euclidean Coxeter groups up to rank $9$, which includes
all of the sporadic examples.  Then special properties of the infinite
families are used to complete the proof.  It would, of course, be more
desirable to have case-free proof of the lattice property, that
project has not yet been completed.

The second main result of \cite{McSu-artin-euclid} establishes that
the crystallographic Garside group $G$ has the structure of an
amalgamated free product, and as a consequence, the natural map from
corresponding euclidean Artin group $A$ to the crystallographic
Garside group $G$ is injective.

\begin{main}[Artin groups are subgroups]\label{main:subgroup}
  For each irreducible euclidean Coxeter group $W=\cox(\wt X_n)$ and
  for each choice of Coxeter element $w$, the crystallographic Garside
  group $G=\gar(\wt X_n,w)$ is an amalgamated free product of explicit
  groups with the euclidean Artin group $A=\art(\wt X_n)$ as one of
  its factors.  In particular, the euclidean Artin group $A$ injects
  into the crystallographic Garside group $G$.
\end{main}

In terms of the groups defined in the previous section the
crystallographic Garside group $G$ is an amalgamated product of $F_w$
and $A$ over $D_w$.  This also means that all of the group
homomorphisms on the top level of Figure~\ref{fig:ten-gps} are
injective.  Injectivity is a consequence of our structural analysis
and not something that was immediately obvious from the definitions of
the maps.  The final result of \cite{McSu-artin-euclid} uses these
embeddings of euclidean Artin groups into crystallographic Garside
groups to elucidate their structure.

\begin{main}[Structure of euclidean Artin groups]
  Every irreducible euclidean Artin group $A=\art(\wt X_n)$ is a
  torsion-free centerless group with a solvable word problem and a
  finite-dimensional classifying space.
\end{main}

Most of these structural results follow immediately from
Theorem~\ref{main:subgroup}.  The only aspect that required a bit more
work is the center.  The Garside structure on $G$, the product
structure on $F_w$, and the fact that we are amalgamating over $D_w$
are all used in the proof that shows the center of $A$ is trivial.
See \cite{McSu-artin-euclid} for details.  

To conclude this survey, I would like to highlight some of the
questions that these results raise.  Now that we understand the word
problem for the euclidean types, can we devise an Artin group
intrinsic solution that avoids the introduction of the
crystallographic Garside groups?  Or perhaps the crystallographic
Garside groups we define are merely the first instance of a natural
geometric completion process?  What about hyperbolic Artin groups and
beyond?  Are there similar procedures that work in these more general
contexts?

\medskip\noindent
\textbf{Acknowledgements:} I would like to thank the organizers of the
2013 Durham symposium on geometric and cohomological group theory for
their invitation to speak about these results and also the organizers
of the 2003 Durham symposium since it was during that earlier
conference that this project was initially conceived.  My
long-suffering collaborators, Noel Brady, John Crisp and Rob Sulway,
also deserve a special note of thanks for putting up with me as I
refused to let this project die a quiet death during its darkest days.

\newcommand{\etalchar}[1]{$^{#1}$}
\def\cprime{$'$}
\providecommand{\bysame}{\leavevmode\hbox to3em{\hrulefill}\thinspace}
\providecommand{\MR}{\relax\ifhmode\unskip\space\fi MR }
\providecommand{\MRhref}[2]{%
  \href{http://www.ams.org/mathscinet-getitem?mr=#1}{#2}
}
\providecommand{\href}[2]{#2}

\end{document}